\documentclass[a4paper,draft]{article}
\usepackage{graphicx}
\usepackage{amsmath,amssymb,amsthm}
\usepackage[textwidth=16cm]{geometry}
\usepackage{todonotes}
\usepackage{enumitem}
\usepackage{hyperref}
\usepackage{bm}
\usepackage{textgreek}

\title{Dimension reduction for a coupled 
electro-elastic\\ saddle-point problem at finite strains}
\author{Kateryna Buryachenko\footnote{%
Humboldt-Universität zu Berlin, 
Institut für Mathematik,
Unter den Linden 6
10099 Berlin, Germany, }%
\and Annegret Glitzky\footnote{%
Weierstrass Institute Berlin, 
Anton-Wilhelm-Amo-Straße 39, 
10117 Berlin, Germany}
\and
Matthias Liero\footnotemark[2] \footnote{Corresponding author: \url{matthias.liero@wias-berlin.de}} 
\and Barbara Zwicknagl\footnotemark[1]}

\date{\today}

\newtheorem{theorem}{Theorem}[section]
\newtheorem{lemma}[theorem]{Lemma}
\newtheorem{proposition}[theorem]{Proposition}
\newtheorem{corollary}[theorem]{Corollary}
\theoremstyle{definition}
\newtheorem{remark}[theorem]{Remark}
\newtheorem{definition}[theorem]{Definition}

\newcommand{\DIV}{\mathrm{div}\,}
\newcommand{\eps}{\varepsilon}
\newcommand{\pl}{\partial}
\newcommand{\calE}{\mathcal{E}}
\newcommand{\calF}{\mathcal{F}}

\newcommand{\calY}{\mathcal{Y}}
\newcommand{\calV}{\mathcal{V}}

\newcommand{\calM}{\mathcal{M}}
\newcommand{\calN}{\mathcal{N}}
\newcommand{\calQ}{\mathcal{Q}}

\newcommand{\calI}{\mathcal{I}}
\newcommand{\ol}{\overline}

\newcommand{\bfkappa}{\boldsymbol{\kappa}}
\newcommand{\bfk}{\boldsymbol{k}}
\newcommand{\dd}{\,\mathrm{d}}
\newcommand{\Cof}{\mathrm{Cof}\,}
\newcommand{\R}{\mathbb{R}}
\newcommand{\N}{\mathbb{N}}
\newcommand{\wt}{\widetilde}
\newcommand{\wh}{\widehat}
\newcommand{\SO}{\mathrm{SO}}
\newcommand{\dist}{\mathrm{dist}}
\newcommand{\curl}{\mathrm{curl}}
\newcommand{\bbK}{\mathbb{K}}

\DeclareMathOperator*{\argmin}{arg\,min}
\DeclareMathOperator*{\esssup}{ess\,sup}

\numberwithin{equation}{section}

\begin{document}

\maketitle
\begin{abstract}
We study the finite deformation of a thin, elastically heterogeneous sheet subject to electrostatic coupling. The interaction between mechanics and electrostatics is formulated as a saddle-point problem involving the deformation and the electrostatic potential. Starting from a three-dimensional electro-elastic model with prestrain in the elastic energy, we rigorously derive a reduced plate model in the bending regime.
To perform the dimension reduction, that is, to derive the energy of a thin object by taking a suitable limit as its thickness tends to zero, we apply $\Gamma$-convergence-type methods to the underlying saddle-point problem. In the case of bivariate functionals, this convergence  is understood  in an adapted epi/hypo-convergence sense. In this concept, we demonstrate the convergence of the rescaled electro-elastic problems to an effective two-dimensional bending model coupled to electric effects. We verify that cluster points of saddle points are saddle points for the limit. 
\end{abstract}



\renewcommand\contentsname{\centering Content}

\tableofcontents
\section{Introduction}

The study of dimension reduction in continuum mechanics provides a 
rigorous framework for deriving effective lower-dimensional theories for thin domains that capture the dominant behavior while reducing e.g.\ computational complexity, see  \cite{ciarlet_destuynder79,  CL1989JBCC, le-dret, M1988SPSS, morgenstern59}. More recently, 
approaches based on variational methods such as $\Gamma$-convergence have proven to be 
a powerful and versatile tool especially for dimension reduction in nonlinear elasticity.
The main task herein, i.e., the convergence of global minimizers of suitably scaled elastic energy functionals, 
has been discussed by many authors since the pioneering works by Le Dret and Raoult \cite{le-dret}, and by Friesecke, James and M\"uller \cite{FJM1,FJM2}.
This problem has gained increasing attention in the case of pre-strained bodies, and several results appeared not only  for the dimension reduction from three- to two-dimensional problems (see e.g.\ \cite{ag-ago19,Bhatt16,Padilla-Garza2022}) but also in the three- to one-dimensional case, i.e., rods 
(see, for instance, \cite{BGKN2022MDNL,BNS2020DHBT, CiRuSo2017GLMP} and references therein).

In this paper, we deal with the derivation of an effective plate model in the bending regime for a coupled three-dimensional electro-elastic model when the 
relative thickness $0<\eps\ll 1$ of the plate goes to zero. 
The model describes the finite deformation of a thin elastically heterogeneous sheet in response to an electric field. While the former is given in terms of the deformation $y:\Omega\to \R^3$ of the reference configuration $\Omega$, the latter is given by the electrostatic potential $\varphi:\Omega\to \R$ that solves the Poisson equation, namely,
\[
-\DIV\big(\bfkappa(x,\nabla y) \nabla\varphi\big) = 
e_0n_\mathrm{ch}(x), 
\]
where $\bfkappa(x,\nabla y)\in\R^{3\times 3}$ is the symmetric and positive definite permittivity tensor in the Lagrangian frame that depends on the deformation gradient, and $n_\mathrm{ch}(x)$ is a fixed charge distribution ($e_0>0$ denotes the elementary charge). In particular, the nonlinear dependence of $\bfkappa$ on $\nabla y$ leads to significant mathematical challenges, see e.g.\ Subsection~\ref{ss34}. We make the crucial scaling assumption, that the electrostatic potential $\varphi$ and the charge density $n_\mathrm{ch}$ are of order $\eps$ (comp.\  \eqref{eq:scalingAssu}).  

Electro-elastic models in the setting of large deformations are highly relevant due to their broad range of applications in describing  electromechanical effects in polymeric materials. Elastomeric materials are sensitive to electric fields and can be used in transducer devices such as actuators and sensors, see e.g.\ \cite{CCMCGCCRFML2023SMMB}. 
This development in materials science requires in turn a development of the mathematical theory to improve the understanding of the electro-mechanical (in particular, electro-elastic) interactions for material characterization and prediction via mathematical analysis and numerical simulations. A thorough study, review of key experiments, discussion of the range of applications, and the history of development of the nonlinear theory of such models have been done, for example, by Dorfmann and Ogden in \cite{DorOgd2017NEEM}. In addition, special attention is given to the development of electro-elastic theory, constitutive equations for electro-elastic materials and then their specialization to isotropic electro-elasticity.  This   is necessary for material characterization and analysis of general electro-elastic coupling problems.
Another novel polyconvex  transversely-isotropic invariant oriented model of Electro-Active Polymers (EAPs) is studied in \cite{HGOK2023PCTI}.
In  that  paper,  a series of numerical examples modelling the performance of transversely isotropic EAPs at large strains are presented, and existence of minimizers, material stability (ellipticity) of  simulations is ensured by means of  a so-called $\cal A$-polyconvexity condition.  Moreover, a   polyconvex basis of invariants for the creation of polyconvex invariant-based constitutive models is also introduced. In the work \cite{MiVaZa2015CSMS}, Miehe, Vallicotti, and Zäh outline variational-based definitions for structural and material stability for EAPs. Herein an enthalpy-based saddle-point principle  is considered as the most convenient setting for numerical implementation. Stability criteria for a canonical energy minimization principle of electro-elasto-statics are formulated and shifted  over to representations related to this enthalpy-based saddle-point principle. 
A linearized version of the EAP model in  \cite{MiVaZa2015CSMS} was studied by  Kru\v{z}\'{i}k and Roub\'{i}\v{c}ek in \cite[Sect. 5.6]{ag-kru19}, see also \cite{Roub2017VMSD}.

In our problem, the main mathematical challenge lies in the fact that it is not formulated as a family of minimization problems, but as a family of saddle-point problems (see, for instance, \cite{AttouchWets83}) 
for the energy functionals $\calF_\varepsilon$. 
This is in contrast, e.g.,   to the setting considered by 
Bartels et al.\ \cite{ag-bar23}. Therein, the authors also deal with a coupled problem
for elastomers, where the deformation $y$ is coupled to the behavior of a director field  via a spontaneous curvature term. In their case, however, the problem is formulated as a joint minimization problem for $y$ and the director field.
In our setting, the energy functional has the form 
\begin{eqnarray}\label{eq:funct}
\calF_\varepsilon(y,\varphi)=\calM_\varepsilon(y)-\calE_\varepsilon (y,\varphi),
\end{eqnarray}
where $\calM_\eps$ denotes the purely elastic part depending only on  $y$, see \eqref{eq:mechEnergyScaled}, while $\calE_\eps$ represents the electrostatic contribution
depending on both, $y$ and the electrostatic potential $\varphi$, see \eqref{eq:elecEnergyScaled}. 

In contrast to linearized piezoelectric problems, see e.g.\ \cite[Sect. 5.6]{ag-kru19}, 
the existence of a saddle point for fixed $\eps$ is not clear as physical principles
require the energy $\calE_\eps$ to be nonconvex (at most polyconvex) 
with respect to the deformation. Thus, classical results such as \cite[Ch.\ 4, Prop.\ 2.2]{EkeTem1999CAVP} cannot be applied. In this regard, the derivation of an effective and to some extent simpler model for thin domains is much desired.

The idea of applying $\Gamma$-convergence methods to saddle-point problems is not new.
A notion of convergence for a sequence of bivariate functionals  $\calF_\eps:\calY\times \calV\to \ol\R:=\R\cup\{\pm\infty\}$ to some limiting functional $\calF_0:\calY_0\times \calV_0\to  \overline{\R}  $  in some metric spaces $\calY,\,\calV,\, \calY_0,\, \calV_0$ is called epi/hypo-convergence and was introduced by Attouch and Wets in 
\cite{AW1981ACNO, AttouchWets83}.  
We build on this idea but suitably adapt it to our setting in order to deal with potential cancellation effects in the energy \eqref{eq:funct} which is defined as difference of two non-negative contributions and a linear part, see Section \ref{sec:dimred}.

Another challenge is that we take into account that the elastic contribution to the total 
energy contains a so-called "prestrain". 
The latter arises, for example, from a layered material structure (see e.g.\ \cite{LCKFBG2018NBBS} for prestrained nanorods) and is common in biological applications \cite{ag-ago19, Lucic2018PhDthesis}.
Typically, these materials are modeled by three-dimensional energy densities 
of the form $\wh W(x,\nabla y(x)) = W_\mathrm{el}(x,\nabla y(x) M(x)^{-1})$ for the deformation gradient $\nabla y$ and a given prestrain $M(x)\in\mathrm{GL}(3)$.
We will build upon the methods in  \cite{ag-ago19,Padilla-Garza2022} to prove the convergence to the bending model for the elastic part.

 The   main result of the present paper is contained in Section~\ref{sec:dimred} (Subsection
~\ref{ss:mainresult}, Theorem ~\ref{thm:main}), and  roughly states the following. Let $(y_\eps^*,\varphi_\eps^*)\in \calY_1\times \calV_1$  be a given sequence of saddle points for the family of functionals of total energy $\calF_\eps:\calY_1\times\calV_1\to \ol\R:=\R\cup \{\pm\infty\}$ (see \ref{eq:abstractSaddleProblemFamily}), such that 
$(y_\eps^*,\varphi_\eps^*)\to (y^*,\varphi^*)\in\calY_0\times\
\calV_0.$ We prove that 
$\calF_\eps$ converges  to $\calF_0:\calY_0\times \calV_0\to \ol \R$ in a sense, which is closely related to epi/hypo-convergence, and which guarantees that $(y^*,\varphi^*)$ is a saddle point of the effective bending model. The specific form of our functionals requires some suitable choices in the topologies and spaces used.

Outline of the paper: In Section~\ref{sec:model}, we introduce the underlying three-dimensional electro-elastic model as well as 
scaling assumptions and non-dimensionalization, see Subsections \ref{ss:3DenergyFunctional} and \ref{ss:scalingNondim}. 
The two-dimensional limit model with limit functional $\calF_0=\calM_0-\calE_0$  is introduced in Subsection~\ref{ss:limitModel}. Moreover, the rigorous assumptions for the limit passage are collected in Subsection \ref{ss24}. 
Section~\ref{sec:analyticRes} provides properties of saddle points for fixed $\varepsilon$. Uniform estimates for saddle points with respect to $\eps$ 
are derived in Subsection \ref{ss:boundedness} (cf.\ Lemma~\ref{lem:apri} and Proposition \ref{prop:conv}).
Section~\ref{sec:dimred} is devoted to the actual dimension reduction. We introduce the  generalized   notion of epi/hypo-convergence for saddle functionals and prove the main result of this paper, Theorem~\ref{thm:main}, on epi/hypo-convergence of the functional $\calF_{\varepsilon}$ of total energy. 
The proof of Theorem~\ref{thm:main} 
relies on establishing several auxiliary asymptotic upper and lower bounds (see Lemmas~\ref{lem:I},
 \ref{lem:IV}, \ref{lem:II}, and \ref{lem:III}). 
 In Subsection \ref{ss34}, we collect some convergence properties that are relevant in finding recovery sequences. 
 Section~\ref{sec:conclusion} is devoted to concluding remarks.

\vskip1cm
\section{Setup for the coupled electro-elastic model}
\label{sec:model}
In this section, we collect all necessary setup statements, introduce the electro-elastic bending model and its limiting model, provide scaling
assumptions and non-dimensionalization for introduced model.

\subsection{Notation}
For the reader's convenience, we collect the notation used throughout the paper here.
{\small
\begin{itemize}
\item
For $x\in\R^3$, we write $x=(x',x_3)$ or $x=(x',t)$, where $x'$ denotes the in-plane components and $x_3$ or $t$ corresponds to the out-of-plane component of $x$;
\item
$\nabla':=(\partial_1,\partial_2)$ means the in-plane gradient and 
$\nabla_\eps:=(\nabla',\frac{1}{\eps}\partial_3)$ denotes the scaled gradient;
\item 
$\{e_1,e_2,e_3\}$ is the standard basis in $\R^3$;
\item
$\R^{n\times n}$ is the vector space of real $n\times n$ matrices, $\mathbb{I}_n\in \R^{n\times n}$ denotes the identity matrix;
\item 
$\curl \,{\bm v}$ stands for the vector product of  $\nabla$ and the vector field ${\bm v}$, i.e., $\curl \,
{\bm v}=\nabla \times {\bm v}$;
\item $\mathrm{GL}(n)$ is the general linear group of degree 
$n$, i.e., the set of 
 invertible $n\times n$ matrices;
\item
$\mathrm{Sym}(n):=\{M\in\R^{n\times n}:\, M^\top=M\}$, the vector space of symmetric matrices,
where $M^\top$ is the transposed matrix of $M$ and $M^{-\top} = (M^{\top})^{-1} = (M^{-1})^\top$;
\item
$\mathrm{Skew}(n):=\{M\in\R^{n\times n}:\, M^\top=-M\},$ the set of skew-symmetric matrices;
\item $M_{\mathrm{sym}}:=\frac{1}{2}(M+M^\top)$, $\nabla_{\mathrm{sym}} y:=\frac{1}{2}(\nabla y+(\nabla y)^\top)$ for vector functions $y$;
\item
$\det M$ is the determinant of $M\in\R^{n\times n}$, and 
$\Cof M = (\det M) M^{-\top}$ denotes the cofactor matrix of $M$;
\item
$\mathrm{SO}(3):=\{M\in \R^{3\times 3}:\, M^\top M=\mathbb{I}_3,\,\mathrm{det} (M)=1\}$ is the set of rotations of $\R^3$;
\item
$W^{2,2}_{\mathrm{iso}}(\omega,\R^3):=\{y\in W^{2,2}(\omega,\R^3):\, \nabla' y^\top \nabla' y=\mathbb{I}_2\}$. For $y\in W^{2,2}_{\mathrm{iso}}(\omega,\R^3)$ let $\nu_y:=\partial_1 y \wedge \partial_1 y$, 
$R_y:=(\partial_1 y|\partial_2 y|\nu_y)$, and 
$\Pi_y:= \nabla' y^\top \nabla' \nu$ the second fundamental form of the surface parametrized by $y$ in local coordinates;
\item 
$L^p_\text{av}(\Omega_1,\R^d)$ stands
for (vector-valued) functions with vanishing average, i.e., $f\in L^p_\text{av}(\Omega_1,\R^d)$ if $f\in L^p(\Omega_1,\R^d)$ and $\int_{\Omega_1}f\dd x =0$. Moreover, we set $W^{k,p}_\text{av}(\Omega_1,\R^d) = W^{k,p}(\Omega_1,\R^d)\cap L^p_\text{av}(\Omega_1,\R^d)$;
\item 
$M_{2\times 2}$ denotes the $2\times 2$ submatrix of 
$M\in \R^{3\times 3}$ resulting from $M$ by omitting the last row and last column;
\item For vectors $a_i,b_i\in \R^n$, their tensor product $a\otimes b\in\R^{n\times n}$ is given by $(a\otimes b)_{ij} = a_i b_j$;
\item
For given $G\in \R^{2\times 2}$, let $G^\circ \in \R^{3\times 3}$ 
denote the matrix obtained by 
\begin{equation} \label{eq:circNotation}
G^\circ:=\left(\begin{array}{cc}
&0\\[-1.5ex]
G&\\[-1.5ex]
&0\\[-0.2ex]
0\quad 0 & 0\\
\end{array}\right);
\end{equation}
\item $c$ stands for a generic constant, $c_\eps$ for a generic constant depending on $\eps$.
\end{itemize}
}

\subsection{Three-dimensional energy functional}
\label{ss:3DenergyFunctional}
We fix a domain $\Omega_h:=\omega\times \left(-h/2,+h/2\right)$, 
where we assume that 
\begin{equation}
\label{eq:domain2D}
\omega\subset\R^2\text{ is an open, bounded and convex domain with piecewise }C^1\text{-boundary}.
\end{equation}
We consider the free energy functional for the deformation $y$ and the electrostatic potential $\varphi$ that consists of a purely mechanical part and an electrostatic contribution, namely
\begin{equation}\label{eq:fullCoupledEnergyH}
\begin{split}
\calF(y,\varphi) &= \calM(y) - \calE(y,\varphi),\quad\text{with}\\
\calM( y) &= \int_{\Omega_h}\Big\{ W_\text{el}(x,\nabla y (x) M(x)^{-1})\det(M(x)) + H(\nabla^2y(x))\Big\}\dd x,\\
\calE(y,\varphi)&=\frac{1}{2}\int_{\Omega_h} 
\Big\{\bfkappa(x,\nabla y) \nabla \varphi(x)\cdot  \nabla \varphi(x)- e_0 n_\mathrm{ch}(x) \varphi(x)\Big\}\dd x,
\end{split}
\end{equation}
where $e_0>0$ is the elementary charge, $n_\mathrm{ch}$
is a fixed charge-density, $M(x)\in\mathrm{GL}(3)$ is the prestrain, and $\bfkappa:\Omega_h\times\mathrm{GL}(3)\to \mathrm{Sym}(3)$
is the symmetric and uniformly positive definite permittivity tensor pulled-back to the reference configuration, namely, for $\bfk\in L^\infty(\Omega_h,\mathrm{Sym}(3))$
uniformly positive definite, we have
\begin{equation}\label{eq:pullbackkappa}
\bfkappa(x,F) = \det(F) F^{-1}\bfk(x)F^{-\top},\quad \text{for almost all }x\in\Omega_h\text{ and for all }F\in\mathrm{GL}(3).
\end{equation}
Note that we have also included a higher order contribution $H(\nabla^2 y)$ to the mechanical energy that acts as a regularization, see Section \ref{sec:analyticRes}.

We are looking for deformations $y$ and electrostatic potentials $\varphi$ such that $y$ minimizes $\calF$ and $\varphi$ is a maximizer in suitable classes of functions,
i.e., we are looking for saddle points of the functional $\calF$. We have the following definition.

\begin{definition}(\cite{AW1981ACNO, AttouchWets83}).\label{def:saddlePoint}
Let $\wt\calY$ and $\wt\calV$ be metric spaces, and $ \calF:\wt\calY\times\wt\calV\to \ol\R$ be a bivariate functional.
A pair $(y_*,\varphi_*)\in \wt\calY\times\wt\calV $ such that  
\begin{equation}\label{eq:saddlePoint}
\calF(y_*,\varphi) \leq \calF(y_*,\varphi_*) \leq \calF(y,\varphi_*)\quad\text{for all }(y,\varphi)\in \wt\calY\times\wt\calV,
\end{equation}
is called a {\it saddle point for the functional $\calF$ on the metric space} $\wt\calY\times\wt\calV$.
\end{definition}

We will use this concept of {\em saddle points} for the functionals $\calF_\eps$ and the limiting functional $\calF_0$ on different pairs of metric spaces $(\wt\calY,\wt\calV)$ and $(\calY_0,\calV_0)$ respectively.

The existence of saddle points for the functional $\calF$ defined in \eqref{eq:fullCoupledEnergyH} is not clear. In particular, standard existence results 
rely on the property that $\calF$ is convex in $y$ and concave in $\varphi$, see \cite{EkeTem1999CAVP,Stefanelli24}, 
while the concavity with respect to $\varphi$ holds, 
the convexity in $y$ is missing due to the 
dependence of $\bfkappa$ on $\nabla y$.
Moreover, the frame indifference principle demands that also the stored-elastic-energy density $W_\text{el}$ is non-convex.
In the case of linearized elasticity, the problem becomes convex-concave, and the existence of saddle points can be shown, \cite[Ch.\ 4, Prop.\ 2.2]{EkeTem1999CAVP}.
We also refer to \cite[Sect. 5.6]{ag-kru19}.

The corresponding Euler--Lagrange equation for the functional $\calF$ with respect to the deformation reads
\[
-\DIV\Big(\underbrace{(\det M)\pl_FW_\text{el}(\nabla yM^{-1})M^{-\top}}_{\text{elastic stress}}\underbrace{-\DIV(\pl_G H(\nabla^2 y))}_\text{hyperstress}+\underbrace{\Sigma_\text{Max}(\nabla y, \nabla\varphi)}_\text{Maxwell stress}
\Big)= 0,
\]
where $\Sigma_\text{Max}$ is the Maxwell stress given for $F=\nabla y$ by
\[
\Sigma_\text{Max}=\big(\bfk F ^{-\top}\nabla \varphi\otimes F^{-\top}\nabla\varphi -\frac{1}{2}\bfk(F^{-\top}\nabla\varphi)\cdot (F^{-\top}\nabla\varphi)\mathbb{I}_3
\big)\Cof F.
\]
The Piola transform $T^E \mapsto T^L:=T^E\Cof \nabla y$ maps an Eulerian tensor field $T^E:y(\Omega_h)\to \R^{3\times 3}$ to its Lagrangian counterpart $T^L:\Omega_h\to\R^{3\times 3}$ 
such that $\det(\nabla y)\DIV_y(T^E) = \DIV_x (T^L) $. Thus, the Maxwell stress takes the more familiar form in spatial (Eulerian) coordinates
\[
\Sigma_\text{Max}^E = -\frac{1}{2}\bfk\nabla_y\phi\cdot \nabla_y\phi \,\mathbb{I}_3
+\bfk \nabla_y\phi\otimes \nabla_y\phi,
\]
where $\phi(y(x)) = \varphi(x)$ denotes the Eulerian electrostatic potential. The Euler--Lagrange equation 
for the electrostatic potential, also called Poisson equation in this case, reads
\[
-\DIV\big((\det F) F^{-1}\bfk(x) F^{-\top} \nabla\varphi\big) = 
e_0n_\mathrm{ch}(x).
\]
For simplicity, we complement  this equation by the condition $\int_{\Omega_h}\varphi\dd x =0$. We refer to \cite[Subsection~2.3]{ag-bar23}, where also additional boundary conditions for $y$ are discussed.

\subsection{Scaling assumptions and non-dimensionalization}\label{ss:scalingNondim}

In order to non-dimensionalize the energy in \eqref{eq:fullCoupledEnergyH}, we introduce a reference length scale $\ell>0$, a reference voltage $V_\text{r}>0$, a reference carrier density $n_\text{r}>0$, and a reference energy density  $E_*>0$ (unit Joule/meter$^3$). We then set the non-dimensional and scaled quantities
\[
\wt\omega:=\ell^{-1}\omega,
\quad \eps = \ell^{-1}h, \quad
\wt x = \ell^{-1} (x',\eps^{-1} x_3) = (\ell^{-1}x',h^{-1}x_3)\in \wt\omega\times\left(-\frac12,+\frac12\right).
\]
The rescaled deformation $\wt y$, electrostatic potential $\wt \varphi$,
and charge density $\wt n_\mathrm{ch}$ are defined via
\[
\begin{gathered}
 \wt y(\wt x) = \ell^{-1} y(\ell \wt x', \ell\eps \wt x_3),\quad \wt \varphi(\wt x) = V_\text{r}^{-1}\varphi(\ell x',\ell \eps \wt x_3),\quad 
\wt n_\mathrm{ch}(\wt x) =n_\text{r}^{-1}{n_\mathrm{ch}(\ell x',\ell\eps\wt x_3)},\\
\wt W_\text{el}(\wt x, F ) = E_*^{-1}W_\text{el}((\ell \wt x',\ell\eps \wt x_3),F),\quad
\wt M(\wt x) = M(\ell \wt x',\ell\eps\wt x_3),
\\
\wt\bfkappa(\wt x,F) = \kappa_0^{-1}\bfkappa (\ell x',\ell\eps\wt x_3,F),\quad \wt H(G) = E_*^{-1} H(G),
\end{gathered}
\]
where $\kappa_0>0$ is the vacuum permittivity.
This scaling leads to the identities
$\nabla_x y = \mathrm{diag}(1,1,\eps^{-1})\nabla_{\wt x}\wt y=:\nabla_\eps\wt y$ and analogously for 
$\varphi$ and $\wt\varphi$. For the Hessians 
of $y$ and $\wt y$, we have the relations
\[
\frac{\partial^2 y}{\partial x_i\partial x_j}  =
\alpha_{ij}^\eps \frac{\partial^2\wt y}{\partial \wt x_i\partial \wt x_j} ,\quad \text{where}
\quad
\alpha_{ij}^\eps= 
\begin{cases}
1/(\ell\eps^2) &\text{if }i=j=3,\\
1/(\ell\eps) & \text{if }(i= 3\text{ and }j\neq 3)\text{ or }(j= 3\text{ and }i\neq 3),\\
1/\ell&\text{otherwise}.
\end{cases}
\]
We will write $\nabla^2_\eps \wt y$, where 
$(\nabla^2_\eps y)_{ijk} = \alpha_{ij}^\eps \frac{\partial^2\wt y_k}{\partial \wt x_i\partial \wt x_j}$, $i,j,k=1,2,3$.

Plugging these identities into the energy functional $\calF$ in \eqref{eq:fullCoupledEnergyH}, we obtain
\[
\frac{1}{\eps^3\ell^3 E_*}\calF(y,\varphi) = \wt\calF_\eps(\wt y,\wt\varphi),
\]
with the rescaled energy functional
\begin{multline}\label{eq:fullCoupledEnergyScaled}
\wt\calF_\eps(\wt y,\wt\varphi) =  \frac{1}{\eps^2}\int_{\Omega_1} 
\wt W_\text{el}(\wt x,\nabla_\eps \wt y (\wt x) \wt M(\wt x)^{-1})\det(\wt M(\wt x)) + \wt H(\nabla_{\eps}^2\wt y(\wt x))\dd\wt x\\
+\frac{e_0V_\text{r}n_\text{r}}{ \eps^2 E_* }\int_{\Omega_1}  
\wt n_\mathrm{ch}(\wt x) \wt \varphi(\wt x)
\dd \wt x
-\frac{\kappa_0V_\text{r}^2}{2\eps^2 E_*}\int_{\Omega_1} 
 \wt\bfkappa\big(\wt x,\nabla_\eps \wt y(\wt x)\big) \nabla_\eps \wt\varphi(\wt x)\cdot \nabla_\eps \wt \varphi(\wt x)\dd \wt x.
\end{multline}

For the reference charge-carrier density $n_\text{r}$ and the reference voltage $V_\text{r}$, we suppose the following smallness assumption
    \begin{equation}
    \label{eq:scalingAssu}
    n_\text{r} = 
    n_\text{r}^\eps = \frac{\eps^2}{\mathcal{L}^3(\Omega_\eps)} = \frac{\eps}{\ell^3}
    \quad\text{and}\quad V_\text{r}=V_\text{r}^\eps = \eps V_*.
\end{equation}
Note that this assumption leads to the factors in front of the last two integrals in \eqref{eq:fullCoupledEnergyScaled} to be of order 1, namely:
        \begin{equation}
        \gamma:=\frac{e_0V_\text{r} n_\text{r}}{\eps^2 E_* } = \frac{e_0V_* n_*}{E_* }\quad\text{and}\quad\beta:=\frac{\kappa_0V_\text{r}^2}{ \eps^2E_*} = \frac{\kappa_0V_*^2}{ E_*}.
        \end{equation}
Concerning the prestrain, we follow \cite{ag-ago19} (and Sect.~4 in \cite{Padilla-Garza2022}) and assume that $\wt M = \mathbb{I}_3+\eps \wt B(x)$, where $\wt B(\wt x)\in\R^{3\times 3}$ is symmetric, see also Assumption~\ref{assu:Prestrain} below.

In the remaining text, we will drop the tilde-notation for notational simplicity.

\subsection{The two-dimensional limit model}\label{ss:limitModel}

In this section, we describe the effective lower dimensional electro-elastic model. As usual in the theory of bending models, see  e.g.\ \cite{ag-ago19,ag-bar23}, let $Q_3(x,\cdot):\R^{3\times 3}\to \R$ be a quadratic form such that
\[
\Big|W_\text{el}(x, \mathbb{I}_3+F)-{\frac{1}{2}}Q_3(x,F)\Big|\leq |F|^2r_W(|F|),\,\forall F\in \mathbb{R}^{3\times 3},~\forall x\in\Omega_1
\]
(see also Assumption \ref{ass:Q3} below).
Moreover, we introduce for a matrix $X\in\R^{2\times 2}$ the quadratic form $Q_2:\Omega_1\times \R^{2\times 2} \to \R$
\[
Q_2(x',t,X):=\min_{z\in \R^3}Q_3\big(x',t,[X^\circ{+}z{\otimes} e_3]\big)
=\min_{z\in \R^3}D^2 W_\mathrm{el}(x',t,\mathbb{I}_3)
\big([X^\circ+z{\otimes} e_3], 
 [X^\circ+z{\otimes} e_3]\big),
\]
where the notation $X^\circ$ is explained in \eqref{eq:circNotation}. 
Now, let $\overline{Q}_2:\omega\times \R^{2\times 2}\to \R$ be given via
\begin{multline*}
\overline{Q}_2(x',X):=  \min_{s\in\R^{2\times 2}} \int_{-1/2}^{1/2}
Q_2(x',t,tX+s-B(x',t)_{2\times 2})\,\dd t\\
=
\min_{s\in\R^{2\times 2}} \int_{-1/2}^{1/2}
\min_{z\in \R^3}D^2 W_\mathrm{el}(x',t,\mathbb{I}_3)
([tX^\circ{+}s^\circ{-}(B_{2\times 2})^\circ {+}z{\otimes} e_3],[tX^\circ{+}s^\circ{-}(B_{2\times 2})^\circ{+}z{\otimes} e_3])\,\dd t.
\end{multline*}
With the assumptions on $W_\text{el}$ below,  $Q_n$, $n=2,3$, are quadratic forms such that
$A\mapsto Q_n(x,A)$ is positive 
semi-definite on $\R^{3\times 3}$ and positive definite on $\mathrm{Sym}(n)$.
Moreover, we have that
$Q_n(x,F)=0$ for all $F\in \mathrm{Skew}(n)$, and 
$F\mapsto Q_n(x,F)$ is strictly convex on $\mathrm{Sym}(n)$.
We refer to \cite[Lemma 2.A.1.]{Lucic2018PhDthesis}.

Additionally, we introduce the unit normal vector $\nu_y$ to the surface belonging to $y\in W^{2,2}(\omega,\R^3)$ with $\nabla' y^\top\nabla'y=\mathbb{I}_2$ and the matrix 
$R_y \in L^2(\omega,\R^{3\times 3})$ via
\begin{equation}\label{eq:mu}
\nu_y(x'):=\frac{\partial_1 y\times \partial_2 y}{\|\partial_1 y\times \partial_2 y\|}
=\partial_1 y\wedge \partial_2 y\in \R^3,\quad
R_y:=(\nabla'y|\nu_y)\quad
\text{such that}\quad R_y^\top R_y=\mathbb{I}_3.
\end{equation}

We work in the following function spaces
\begin{equation}\label{eq:spacesY0V0}
\begin{split}    
    \calY_0&:=W^{2,2}_\mathrm{iso}(\omega,\R^3):=\big\{y\in W^{2,2}(\omega,\R^3):
\nabla' y^\top \nabla' y=\mathbb{I}_2\big\},\\
\calV_0&:=W^{1,2}_\text{av}(\omega)=\big\{\varphi\in W^{1,2}(\omega):\, \textstyle\int_{\omega}\varphi\,\dd x'=0\big\}.
\end{split}
\end{equation}
By the canonical extension of functions defined on $\omega$ to functions on $\Omega_1$, we can view the spaces $\calV_0$ and $\calY_0$
as subsets of $W^{1,2}(\Omega_1)$ and $W^{2,2}(\Omega_1,\R^{3})$ respectively.

We define for $y\in\calY_0$ the limiting mechanical energy functional  $\calM_0(y)$ by
\begin{equation*}
\begin{split}
\calM_0(y)&:=\frac{1}{2}\int_{\omega}\overline{Q}_2(x',\nabla y^\top\nabla \nu_y)\,\dd x' \\
&=
\frac{1}{2}\int_{\omega}\min_{s\in\R^{2\times 2}}\int_{-1/2}^{1/2}
\min_{z\in \R^3}Q_3(x',t,t\nabla y^\top\nabla \nu_y-(B_{2\times 2})^\circ+s+z\otimes e_3)\,\dd t\,\dd x'.
\end{split}
\end{equation*}

Concerning the limiting electrostatic energy, the effective energy density is obtained
by minimizing the effect of the vertical derivatives, i.e., for fixed $y\in \calY_0$ and associated $R_y= (\nabla' y|\nu_y)$ we define
\begin{equation}\label{eq:defMYPHI}
\begin{split}
m_{y,\varphi}(x')&=\argmin_{\widehat m\in \R}
\bbK_y(x') 
\binom{\nabla'\varphi( x')}{\wh m}\cdot 
\binom{\nabla' \varphi( x')}{\wh m},\text{ where}\\
\bbK_y(x')&:=\int_{-1/2}^{1/2}\bfkappa(x',t,R_y(x'))\dd t =R_y(x')^\top\ol\bfk(x')R_y(x'), 
\\&\qquad\text{with } \ol\bfk(x') := \int_{-1/2}^{1/2}\bfk(x',t)\dd t,
\end{split}
\end{equation}
$R_y$ by construction does not depend on $t=x_3$.
Measurability of $m_{y,\varphi}$ follows from general results for optimal values of normal integrands, see \cite[Theorem 14.37]{RocWet1998VA}.
Furthermore, by the assumptions on $\bfk$ in \ref{eq:AssuElectric2}, $\bbK_y$ is uniformly positive definite and bounded, 
and hence for almost every $x'\in \omega$, we have
\begin{equation}\label{eq:KElliptic}
\begin{split}
c|m_{y,\varphi}(x')|^2&\leq c\big|\big((\nabla'\varphi( x'), m_{y,\varphi}( x')\big)\big|^2
\leq \bbK_y(x')(\nabla'\varphi( x'),
m_{y,\varphi}( x'))\cdot (\nabla' \varphi( x'), m_{y,\varphi}( x'))\\
&\le \bbK_y(x')(\nabla'\varphi( x'),0)
\cdot (\nabla' \varphi( x'),0)\leq \|\bfk\|_{L^\infty(\Omega_1, 
\R^{3\times 3})}|\nabla'\varphi(x')|^2.
\end{split}
\end{equation}
We can explicitly compute $m_{y,\varphi}$ using Schur complements. Let us write
\[
\bbK_y = 
\begin{pmatrix}
\ol \bbK_y & K_y\\
K_y^\top & k_y
\end{pmatrix}
\quad\text{with}\quad \ol \bbK_y(x')\in\R^{2\times 2},~K_y(x')\in\R^2,~k_y(x')\in\R_+.
\]
Then, the minimizer in \eqref{eq:defMYPHI} 
satisfies $m_{y,\varphi}= -\frac{1}{k_y} K_y\cdot\nabla'\varphi$,
and we obtain
\begin{equation}
\begin{split}
&\bbK_y(x') 
\binom{\nabla'\varphi( x')}{m_{y,\varphi}(x')}\cdot 
\binom{\nabla' \varphi( x')}{m_{y,\varphi}(x')}
= \bbK_y^\mathrm{eff}(x')\nabla'\varphi( x')\cdot \nabla'\varphi( x')\\&\quad\text{with}\quad
\bbK_y^\mathrm{eff}(x') 
= \ol \bbK_y(x') - \frac{1} {k_y(x')}K_y(x')\otimes K_y(x')\in \R^{2\times 2}.
\end{split}
\end{equation}
In particular, the effective tensor  $\mathbb{K}^\mathrm{eff}_y$ is 
the Schur complement of the component $k_y$ in $\bbK_y$. 
Obviously, $\bbK_y^\mathrm{eff}$ is symmetric. It is also uniformly
positive definite and bounded, see \eqref{eq:KElliptic}.

For $y\in\calY_0$  and $\varphi\in\calV_0$, the effective 
electrostatic energy functional $\calE_0(y,\varphi) $ is defined  by 
\begin{equation}  
\label{eq:limitElecEnergy}
    \calE_0(y,\varphi) 
    := \frac{\beta}{2}\int_{\omega} 
    \mathbb{K}^\mathrm{eff}_y(x')\nabla'\varphi( x')\cdot \nabla' \varphi( x')\dd  x'
    -\gamma\int_{\omega}  
    \ol n_\mathrm{ch}( x')  \varphi( x') \dd  x',
\end{equation}
where $\ol n_\mathrm{ch}(x')
:= \int_{-1/2}^{1/2}n_\mathrm{ch}(x',t)\dd t$.
Thus,  minimizers for $\calE_0(y,\cdot)$ satisfy the effective Poisson equation
\[
-\DIV\big(\mathbb{K}^\mathrm{eff}_y(x') \nabla' \varphi\big)=\frac{\gamma}{\beta}\ol n_\mathrm{ch}(x')\quad\text{in }\omega.
\]

\begin{remark}[Isotropic case]
In the isotropic case $\bfk=k_*\mathbb{I}_3$ with $k_*>0$, we obtain  $\bbK_y= k_*\mathbb{I}_3$, since $R_y\in\SO(3)$, such that $m_{y,\varphi}\equiv0$. In particular, in this case we also have that $\bbK^\mathrm{eff}_y = k_*\mathbb{I}_2$. The limiting electrostatic energy does not depend on the deformation $y$ in this case and the two equations for the deformation and electrostatic potential decouple.
\end{remark}

Finally, we introduce the effective total free energy $\calF_0$ by
\begin{equation}
    \label{eq:totaleffectiveEnergy}
    \calF_0(y,\varphi) = \calM_0(y) -\calE_0(y,\varphi)\quad \text{for }(y,\varphi)\in \calY_0\times\calV_0.
\end{equation}

\subsection{Assumptions for the bending model}\label{ss24}

We impose the following assumptions that allow us to pass to the bending model.
\begin{enumerate}[label=(A\arabic*)]

    \item \label{assu:W}
    The elastic stored-energy density $W_\text{el}:\Omega_1\times \R^{3\times 3}\to [0,\infty]$ satisfies the conditions: 
\begin{enumerate}[label=(W\arabic*)]
\item \label{assu:inva}
$W_\text{el}(x,RF)=W_\text{el}(x,F),\,\forall F\in \R^{3\times 3}, R\in \SO(3)$ for a.e.\ $x\in\Omega_1$ (frame indifference);
\item \label{assu:WgrowthSO3}
For all $F\in\R^{3\times 3} $ and for a.e.\ $x\in\Omega_1$ it holds that
(non-degeneracy and natural state)
\begin{align}
    W_\text{el}(x,F)&\geq\frac{1}{C_W}{\rm dist^2}(F,\,\SO(3)),\\
    W_\text{el}(x,F)&\leq C_W {\rm dist^2}(F,\,\SO(3))\quad \text{if}\quad {\rm dist^2}(F,\,\SO(3))\leq \frac{1}{C_W};
\end{align}

\item \label{ass:Q3}
For almost every $x\in\Omega_1$, there exists a quadratic form $Q_3(x,\cdot ):\mathbb{R}^{3\times 3}\to \mathbb{R}$ 
such that
\[ 
\Big|W_\text{el}(x, \mathbb{I}_3+F)-\frac{1}{2}Q_3(x,F)\Big|\leq |F|^2r_W(|F|),\,\forall F\in \mathbb{R}^{3\times 3},
\]
where $r_W:[0,\infty)\to[0,\infty]$ is monotone with $\lim_{t\to 0}r_W(t)= 0$;
\item \label{assu:Wgrowth}
There exist $q_W>6$ and $C_W>0$ such that for all $F\in\R^{3\times 3} $ and for a.e.\ $x\in\Omega_1$ there  holds the lower bound
\[
W_\text{el}(x, F)\geq \begin{cases}
  \frac{1}{C_W} \max\big\{|F|^{q_W},\,\det (F)^{-\frac{q_W}{2}}\big\}-C_W  & \text{if }\det (F)>0, \\
  \infty   & \text{else;}
\end{cases}
\]
\item \label{assu:conti}
$W_\text{el}$ is continuous and there exists a neighborhood $\mathcal{U}$ of $\SO(3)$ such that $W_\text{el}\in C^2(\mathcal{U})$ and $D^2W_\text{el}$ is uniformly equicontinuous, i.e., 
\[
\forall \varepsilon>0\,\,\exists \delta>0\,\,\forall\,F_1, F_2\in \mathcal{U} :~ \Big[ |F_1{-}F_2|<\delta ~ \Longrightarrow~|D^2 W_\text{el}(x,F_1)-D^2W_\text{el}(x,F_2)|<\varepsilon\Big].
\]
\end{enumerate}

\item 
    \label{assu:hyper}
    The hyperstress potential $H:\mathbb{R}^{3\times 3\times 3}\to \mathbb{R}_+$ has the form 
    \begin{equation}
    H(G) = H_\eps(G) = \eps^{\alpha_H} H_*(G),
    \end{equation}
    for some exponent $\alpha_H>0$ to be fixed later. The function $H_*$
    is convex, 
    continuously differentiable and there exist $K_H\geq c_H>0$ and $q_H>3$ such that
    $q_W/2 > 3 q_H/(q_H{-}3) > 3$ (for $q_W$ from \ref{assu:Wgrowth}) and
\[
     c_H|G|^{q_H}\le H_*(G)\le K_H(1+|G|^{q_H}) 
     \text{ for all }G\in
    \mathbb{R}^{3\times 3\times 3}.
    \]
    Moreover, $H$ is frame indifferent, meaning that for all $ R\in \SO(3)$
    and $G\in \mathbb{R}^{3\times 3\times 3} $ we have
    $H(R G)=H(G)$.
\item \label{assu:Prestrain}
    For the prestrain $M$, we assume that $
    M(x) = M_\eps(x) :=
    \mathbb{I}_3 + \eps B(x)$ with $B\in L^\infty(\Omega_1,\mathbb{R}^{3\times 3}_\mathrm{sym})$.
\item\label{eq:AssuElectric} 
    For the reference charge-carrier density $n_\text{r}$ and the reference voltage $V_\text{r}$, we suppose the following smallness assumption
    \begin{equation}\label{eq:AssuCharge}
    n_\text{r} = 
    n_\text{r}^\eps = \frac{\eps^2}{\mathcal{L}^3(\Omega_\eps)} = \frac{\eps}{\ell^3}
    \quad\text{and}\quad V_\text{r}=V_\text{r}^\eps = \eps V_*.
\end{equation}

\item\label{eq:AssuElectric2}
The quantity $n_\mathrm{ch}\in L^\infty(\Omega_1)$ is a fixed charge-density, and $\bfkappa:\Omega_1\times\mathrm{GL}(3)\to \mathrm{Sym}(3)$ is defined as in \eqref{eq:pullbackkappa} with $\bfk\in L^{\infty}(\Omega_1,\R^{3\times 3})$ being symmetric and uniformly positive definite such that there exists $ \kappa_0>0 $ 
with $\bfk(x)\xi\cdot\xi\geq \kappa_0|\xi|^2$ for all $\xi\in\R^3$ and almost every $x\in\Omega_1$.
\end{enumerate}

\begin{remark}\label{rem:assu}
    \begin{enumerate}
        \item The conditions \ref{assu:inva}--\ref{ass:Q3} are standard assumptions in the context of the derivation of plate theories from 3d-nonlinear elasticity. We also use Assumptions \ref{assu:Wgrowth} and \ref{assu:conti} as in similar models 
        on the elastic stored-energy density in  \cite{ag-bar23} and \cite{Padilla-Garza2022}.
        However, note that we need the stronger condition $q_W>6$ (see Lemma~\ref{lem:qw6}) than in \cite{ag-bar23}, where $q_W>4$ is sufficient.
        
        \item A typical choice for the hyperstress potential is $H(G)=\frac{\eps^{\alpha_H}}{{q_H}}|G|^{q_H}$.
        In the limit $\eps\to0$, we want its contribution to the free energy to vanish. This means, that the exponent $\alpha_H$ should satisfy $\alpha_H>2+2{q_H}$.

        \item Assumption \ref{assu:Prestrain} for the prestrain $M_\eps$ gives for small $\eps>0$ the expansions 
        \begin{align}
                 M_\eps( x)^{-1} &= \mathbb{I}_3- \eps  B( x) + O(\eps^2),\text{\ and\ }\\ 
                 \label{eq:expansionDeterminant}
                \det  M_\eps(x) &= 1+\eps\,\mathrm{tr}( B( x)) + O(\eps^2).
         \end{align}
         We highlight that the factor $\det M_\eps$ does not usually appear in the literature, see e.g.\ \cite{ag-ago19, Padilla-Garza2022}. Here, we follow the discussion in \cite{BNS2020DHBT} where it is assumed that the prestrained body consists of two materials occupying subdomains $\Omega^{(1)}$ and $\Omega^{(2)}$ with $\ol\Omega = \ol \Omega^{(1)}\dot\cup\,\ol \Omega^{(2)}$. It is assumed that subsets of $\Omega^{(1)}$ or $\Omega^{(2)}$ relax to a stress-free
(elastic energy minimizing) state described by affine deformations $x\mapsto M_i x$, $i=1,2$, such that $\wt \Omega^{(i)} = M_i\Omega^{(i)}$ define stress-free reference configurations for the material. The elastic energy of a deformation $\wt u$ defined relative to $\wt \Omega^{(i)}$ is given by $ \int_{\wt \Omega^{(i)}}W_\text{el}(\nabla \wt y) \dd x$. The multiplicative decomposition
$\nabla y = F_\text{el} M_i$ and the factor $\det M_i$
arise from the change-of-variables formula when going back to the common reference configuration $\Omega$.
A typical example is a crystalline heterostructure consisting of two materials on top of each other with different lattice constants, see e.g.\ \cite{HMLF2024EFM3}. However, note that in our setting  the factor satisfies $\det M_i =1 + O(\eps)$ due to \eqref{eq:expansionDeterminant}.
               
    \end{enumerate}
\end{remark}

The assumptions above lead to the following definition of (subsets of) function spaces for the mechanical deformations and the electrostatic potentials
\begin{equation}\label{eq:spacesY1V1}
\begin{split}
\calY_1 
&:= \big\{ y\in W^{2,q_H}(\Omega_1,\R^3)\,\big|\, (\det \nabla y)^{-1} \in L^{q_W/2}(\Omega_1),~\textstyle \int_{\Omega_1}y\dd x = 0 \big\}, \text{\quad and}\\
\calV_1 &:= \{\varphi\in W^{1,2}(\Omega_1)\,|\, \textstyle\int_{\Omega_1}\varphi \dd x =0\}=W_{{\rm av}}^{1,2}(\Omega_1).
\end{split}
\end{equation}

The purely mechanical energy functional $\calM_\eps:\calY_{ 1}\to [0,\infty]$ is defined as
\begin{equation}
    \label{eq:mechEnergyScaled}
    \calM_\eps(y) := \frac{1}{\eps^2}\int_{\Omega_1} 
    \Big\{
    W_\text{el}( x,\nabla_\eps  y ( x)  M_\eps( x)^{-1})\det( M_\eps( x)) + \eps^{\alpha_H} H_*(\nabla_{\eps}^2 y( x))\Big\}\dd x.
\end{equation}
The electrostatic energy $\calE_\eps:\calY_1\times\calV_1\to \R$ is given by
\begin{equation}\label{eq:elecEnergyScaled}
\begin{split}
    \calE_\eps(y,\varphi) &:= \frac{\beta}{2}\int_{\Omega_1} 
    \bfkappa\big( x,\nabla_\eps  y( x)\big)\nabla_\eps \varphi( x))\cdot \nabla_\eps  \varphi( x)\dd  x-\gamma\int_{\Omega_1}  
    n_\mathrm{ch}( x)  \varphi( x) \dd  x.
\end{split}
\end{equation}
Thus, the total free energy $\calF_\eps:  \calY_1\times \calV_1 \to \R_\infty:=\R\cup\{+\infty\}$ reads
\begin{equation}
    \label{eq:totalEnergyScaled}
    \calF_\eps(y,\varphi) = \calM_\eps(y) -\calE_\eps(y,\varphi)\quad \text{for }(y,\varphi)\in\calY_1\times\calV_1.
\end{equation}

\begin{remark}
For fixed $\eps>0$ and deformations $y\in \calY_1$ 
with finite mechanical energy, i.e., $\calM_\eps(y)\leq C_\mathrm{M} < \infty$, we get from Assumptions \ref{assu:Wgrowth} and \ref{assu:hyper} by the Healey--Krömer theorem \ref{thm:Healey} (see also \cite[Subsec.\ 3.1]{ag-mie20} and \cite[Theorem 2.5.3]{ag-kru19})  that there exists a constant $C_\mathrm{HK}^\eps=C_\mathrm{HK}(C_\mathrm{M},q_H,q_W,\eps)>0$ such that
        \begin{align}
        \begin{gathered}
        \|y\|_{W^{2,q_H}}\leq C_\text{HK}^\eps,\quad \|y\|_{C^{1-3/q_H}}\leq C_\text{HK}^\eps,\quad \|(\nabla y)^{-1}\|_{C^{1-3/q_H}}\leq C_\text{HK}^\eps,\\ 
        \text{and\ }\det\nabla y(x)\geq 1/C_\text{HK}^\eps\quad\text{for all }x\in\Omega_1.
        \end{gathered}
        \end{align}
        Note that in general $C_\text{HK}^\eps\to \infty$ for $\eps\to 0$.

        For fixed $\eps>0$ and for a fixed $y\in\calY_1$ with bounded mechanical energy, the permittivity tensor satisfies $\bfkappa(x,\nabla_\eps y) \in L^\infty(\Omega_1,\R^{3\times 3})$ and is uniformly positive definite,
        i.e., $\bfkappa(x,\nabla_\eps y)\xi\cdot \xi\geq c_\eps|\xi|^2$ for some constant $c_\eps>0$ and all $\xi\in\R^3$. In particular, the electrostatic energy $\calE_\eps$ is finite for every pair $(y,\varphi)\in\calY_1\times \calV_1$. 
\end{remark}

\section{Preliminary analytical results}
\label{sec:analyticRes}

Here, we collect results that will be used in the proof of the limit passage $\eps\to0$  in Section~\ref{sec:dimred}.

\subsection{Properties for fixed thickness}
\label{ss25}

For the dimension reduction in Section~\ref{sec:dimred}, we will assume that  a saddle point $(y_\eps,\varphi_\eps)\in \calY_1\times\calV_1$ for the functional $\calF_\eps$ defined in \eqref{eq:fullCoupledEnergyScaled} exists. In particular, such a pair satisfies
\begin{equation*}
\forall (\wh y,\wh \varphi)\in\calY_1\times \calV_1:\quad\calF_\eps(y_\eps,\wh \varphi) \leq \calF_\eps(y_\eps,\varphi_\eps)
\leq \calF_\eps(\wh y,\varphi_\eps).
\end{equation*}
However, the existence of such a pair is not guaranteed.
 Standard existence results for saddle points typically require, besides compactness properties,
that the functional $(y,\varphi)\mapsto \calF_\eps(y,\varphi)$
is convex in the first and concave in the second variable, see 
e.g.~\cite[Ch.\ 4, Prop.\ 2.2]{EkeTem1999CAVP}. In \cite{Stefanelli24}, the existence of saddle points under weaker coercivity assumptions is established, yet the analysis there also relies on convex–concave structures.
While we have concavity of $\varphi\mapsto\calF_\eps(y,\varphi)$, the convexity of $y\mapsto \calF_\eps(y,\varphi)$ cannot be assumed 
due to the dependence of $\bfkappa$ on $\nabla y$ and since it violates the physical principle of frame invariance.

For fixed $\varepsilon>0$  and any  deformation $y\in\calY_1$ with finite mechanical energy $\calM_\eps(y)<\infty$, 
we prove next the existence of a unique minimizer of the electrostatic energy $\calE_\eps(y,\cdot)$.

\begin{lemma}\label{lem:ExistPoisson}
We assume \ref{assu:W} -- \ref{eq:AssuElectric2}.
Let $\eps>0$ be fixed and consider a deformation $y\in \calY_1$
 with finite mechanical energy, i.e., $\calM_\eps(y)< \infty$. Then, there is a unique weak solution $\varphi=\varphi(y,\eps)\in \calV_1$ to
\begin{equation}
\label{eq:weakPGeps} 
\beta \int_{\Omega_1} \bfkappa (\nabla_\eps y)\nabla_\eps \varphi\cdot \nabla_\eps\ol \varphi\dd x
=
\gamma \int_{\Omega_1}  
n_\mathrm{ch}\ol \varphi \dd x\quad \forall \,\ol\varphi \in \calV_1. 
\end{equation}
Moreover, this $\varphi$ is the unique minimizer of the electrostatic 
energy $\calE_\eps(y,\cdot):\calV_1\to \R$, defined in \eqref{eq:elecEnergyScaled}.

\end{lemma}
 
\begin{proof} 
Since $\calM_\eps(y)< \infty$, the properties of the hyperstress potential in \ref {assu:hyper} and the coercivity  property of the elastic stored energy \ref{assu:Wgrowth} 
ensure that
\[
\nabla_\eps y\in L^{q_W}(\Omega_1,\R^{3\times 3}),\quad
(\det \nabla_\eps y)^{-1}\in L^{q_W/2}(\Omega_1),\quad
\nabla_\eps^2 y\in L^{q_H}(\Omega_1,\R^{3\times 3\times 3}),
\] 
and we find a constant $c_\eps>0$ depending on $\eps$ such that
\[
\|\nabla y\|_{L^{q_W}(\Omega_1,\R^{3\times 3})} +
\|(\det \nabla y)^{-1}\|_{L^{q_W/2}(\Omega_1)}+
\|\nabla^2 y\|_{L^{q_H}(\Omega_1,\R^{3\times 3\times 3})}
\le c_\eps.
\]
The Healey--Kr\"omer theorem \ref{thm:Healey}   ensures the existence of a constant  $C_\text{HK}^\eps>0$ such that
\begin{align*}
        \begin{gathered}
        \| y\|_{W^{2,q_H}}\leq C_\text{HK}^\eps,\quad \| y\|_{C^{1-3/q_H}}\leq C_\text{HK}^\eps,\quad \|(\nabla  y)^{-1}\|_{C^{1-3/q_H}}\leq C_\text{HK}^\eps,\\ 
     \text{and\ }   \det\nabla  y(x)\geq 1/C_\text{HK}^\eps\quad\text{for all }x\in\Omega_1.
        \end{gathered}
        \end{align*}
Therefore, since $q_H>3$, we have $\nabla y\in L^\infty(\Omega_1,\R^{3\times 3})$,
$\det \nabla y\in L^\infty(\Omega_1)$, $\Cof\nabla y\in L^\infty(\Omega_1,\R^{3\times 3})$, and for fixed $\eps>0$ also
$\nabla_\eps y\in L^\infty(\Omega_1,\R^{3\times 3})$ ,
$\det \nabla_\eps y\in L^\infty(\Omega_1)$, and $\Cof\nabla_\eps y\in L^\infty(\Omega_1,\R^{3\times 3})$. Moreover, the Healey--Kr\"omer theorem also yields $\det\nabla_\eps y(x)\geq c_\eps$.

We set $S_\eps := \mathrm{diag}(1,1,1/\eps)$ such that $\nabla_\eps\varphi = S_\eps \nabla \varphi$. Then, 
with Assumption \ref{eq:AssuElectric2} and 
$\Cof\nabla_\eps y=(\nabla_\eps y)^{-1}\det \nabla_\eps y$,
 we can bound the coefficient
 matrix by 
\begin{equation}\label{eq:weak1}
\begin{split}
& \big\|\det \nabla_\eps y \, S_\eps(\nabla_\eps y)^{-1}
\bfk 
(\nabla_\eps y)^{-\top}
S_\eps\big\|_{L^\infty(\Omega_1,\R^{3\times 3})}\\
& =\Big\|\frac{S_\eps\,(\Cof\nabla_\eps y)^\top
\bfk \, \Cof\nabla_\eps y\,S_\eps}
{\det \nabla_\eps y}\Big\|_{L^\infty(\Omega_1,\R^{3\times 3})}
\le c_\eps.
\end{split}
\end{equation}
For $\zeta\in\R^3$ and $F\in\R^{3\times 3}$, we have 
$|\zeta|=|FF^{-1}\zeta|\le |F||F^{-1}\zeta|$. Thus, by \ref{eq:AssuElectric2}
we obtain 
for  all $\zeta\in\R^3$  
\begin{equation}\label{eq:weak2}
\begin{split}
&(\det \nabla_\eps y) \big[S_\eps(\nabla_\eps y)^{-1}
\, \bfk \,
(\nabla_\eps y)^{-\top}S_\eps\zeta\big]\cdot
\zeta\\
&= (\det \nabla_\eps y )\, \bfk\,
\big[(\nabla_\eps y)^{-\top}S_\eps\zeta\big]\cdot
\big[(\nabla_\eps y)^{-\top}S_\eps\zeta\big]\\
&\ge \kappa_0 (\det \nabla_\eps y)
\big|(\nabla_\eps y)^{-\top}S_\eps\zeta\big|^2\\
&\ge 
\frac{\kappa_0 \det \nabla_\eps y}{\|S_\eps^{-1}(\nabla_\eps y)^\top\|_{L^\infty(\Omega_1,\R^{3\times 3})}^2}|\zeta|^2
\ge c_\eps |\zeta|^2\quad \text{ a.e.\ in }\Omega_1.
\end{split}
\end{equation}
In view of  \eqref{eq:weak1} and \eqref{eq:weak2}, the classical elliptic theory (Lax--Milgram lemma) 
ensures the existence of a unique weak solution $\varphi=\varphi(y,\eps)\in \calV_1$ to 
\begin{equation}\label{eq:PGeps}
-\DIV[\det \nabla_\eps y \, S_\eps(\nabla_\eps y)^{-1}
\bfk 
(\nabla_\eps y)^{-\top}
D_\eps\nabla \varphi]=\frac{\beta}{\gamma}n_\mathrm{ch},
\end{equation}
meaning that \eqref{eq:weakPGeps} is fulfilled, which is nothing but the Euler--Lagrange equation for $\varphi\mapsto \calE_\eps(y,\varphi)$.
In other words, $\varphi=\varphi(y,\eps)\in \calV_1$ is the unique minimizer
of the (strictly convex) functional $\varphi\mapsto \calE_\eps(y,\varphi)$ for the fixed $y\in\calY_1$, thus 
\[
\calE_\eps(y,\wh \varphi)\ge \calE_\eps(y,\varphi)
\quad\text{ and }\quad
\calF_\eps(y,\wh \varphi)\le \calF_\eps(y,\varphi)\quad \forall 
\wh \varphi \in \calV_1.
\]
This finishes the proof.
\end{proof}

\begin{corollary}\label{cor:FE}
We assume \ref{assu:W} -- \ref{eq:AssuElectric2}.
Let $\eps>0$ be fixed and $y\in \calY_1$ a deformation. Further, let 
the pair $(y,\varphi)\in  
\calY_1 \times \calV_1$ 
 fulfill
 \begin{equation}\label{eq:le1}
 \calM_\eps(y)<\infty\quad\text{and} \quad
 \calF_\eps(y,\bar\varphi)\le 
 \calF_\eps(y,\varphi)\quad\forall 
 \bar\varphi\in  \calV_1, 
 \end{equation}
i.e., $\varphi$ is a maximizer of $\calF_\eps(y,\cdot)$
(minimizer of $\calE_\eps(y,\cdot)$) for given deformation $y$.
Then 
\begin{equation}\label{eq:PG0}
\beta \int_{\Omega_1} \bfk \big(\nabla_\eps y^{-\top}\nabla_\eps \varphi\big)\cdot \big(\nabla_\eps y^{-\top}\nabla_\eps\varphi\big)\det \nabla_\eps y\dd x
=
\gamma \int_{\Omega_1}  
n_\mathrm{ch}\varphi \dd x, 
\end{equation}  
and the free energy can be rewritten in the two alternatives
\begin{equation}
\label{eq:formFunctionalPoisson}
\begin{split}
\calF_\eps(y,\varphi)
& = \calM_\eps(y) + \frac{\beta}{2} \int_{\Omega_1} \bfk 
\big(\nabla_\eps y^{-\top}\nabla_\eps \varphi\big)\cdot \big(\nabla_\eps y^{-\top}\nabla_\eps\varphi\big)\det \nabla_\eps y
\dd x\\
&
= \calM_\eps(y) + \frac{\gamma}{2} \int_{\Omega_1} 
n_\mathrm{ch}\varphi \dd x.
\end{split}
\end{equation}
\end{corollary}

\begin{proof}
Note that the unique minimizer $\varphi\in \calV_1$ of the (strictly convex) functional $\varphi\mapsto \calE_\eps(y,\varphi)$ for the fixed $y\in\calY_1$ is the weak solution to the Euler--Lagrange equation \eqref{eq:PGeps},
meaning that \eqref{eq:weakPGeps} is fulfilled.
Since $\varphi=\varphi(y,\eps)\in \calV_1$ is an admissible test function in  \eqref{eq:weakPGeps}, we obtain \eqref{eq:PG0}.
The both alternative expressions for the free energy $\calF_\eps(y,\varphi)$ then follow from     \eqref{eq:PG0},
\eqref{eq:elecEnergyScaled}, and \eqref{eq:totalEnergyScaled}.
\end{proof}

The existence of a unique electrostatic potential $\varphi = \varphi(y,\eps)$ that minimizes the electrostatic energy $\calE_\eps$ for a fixed deformation $y\in\calY_1$ with finite elastic energy follows from Lemma~\ref{lem:ExistPoisson}. In particular,
we can always find pairs
$(y_\eps,\varphi_\eps)\in\calY_1\times\calV_1$ such that
\begin{equation}
\label{ass:AlmostSaddlePoint}
\varphi_\eps\text{ maximizes } \wh \varphi\mapsto \calF_\eps(y_\eps,\wh \varphi)\text{ over }\calV_1\quad\text{and}\quad \sup_{\eps>0}\calF_\eps(y_\eps,\varphi_\eps)< \infty.
\end{equation}
In contrast, the existence of a deformation $y_\eps$ that minimizes the total 
energy $\calF_\eps(\cdot, \varphi) = \calM_\eps - \calE_\eps(\cdot,\varphi)$ for a fixed electrostatic potential $\varphi\in \calV_1$ is not trivial, since
$\calE_\eps(y,\varphi)=\infty$ can occur if $\det\nabla y$ is not uniformly bounded away from $0$.

\subsection{Boundedness and converging subsequences}\label{ss:boundedness}

We adopt several ideas from \cite{ag-bar23} for the passage to the limit $\eps\to 0$, i.e., the dimension reduction. 
In particular, Lemma~\ref{lem:apri} below is the analogue of  \cite[Lemma~4.1]{ag-bar23}.
 
\begin{lemma}\label{lem:apri}
We assume \ref{assu:W} -- \ref{eq:AssuElectric2}.
Let $(y_\eps,\varphi_\eps)_{\eps>0}\subseteq \calY_1 \times \calV_1$ 
be a sequence such that 
 \begin{equation*}
 \calF_\eps(y_\eps,\bar\varphi)\le 
 \calF_\eps(y_\eps,\varphi_\eps)\quad\forall 
 \bar\varphi\in   
 \calV_1 \quad \text{ and }\quad
 \calF_\eps(y_\eps,\varphi_\eps)\leq C.
 \end{equation*}
 Then, there exist constants $c>0$ and  $\eps_0>0$ such that for all $\eps \in (0,\eps_0)$
\begin{align}
&\int_{\Omega_1}
\dist^2\big(\nabla_\eps y_\eps M_\eps^{-1}, \SO(3)\big)\det M_\eps\dd x\le C_W\eps^2,\quad
\int_{\Omega_1}\dist^2(\nabla_\eps y_\eps, \SO(3))\,\dd x\le c\eps^2,
\label{eq:conv1}\\
&\int_{\Omega_1}\Big(|\nabla_\eps y_\eps|^{q_W}
+|\det (\nabla_\eps y_\eps)|^{-\frac{q_W}{2}}\Big)\dd x
\le c(q_W,C_W), \label{eq:conv2}\\
&\int_{\Omega_1} |F_\eps^{-\top}\nabla_\eps \varphi_\eps|^2 \det F_\eps\dd x\le c,\quad \text{ where } F_\eps=\nabla_\eps y_\eps,
\label{eq:conv3}\\ 
&\int_{\Omega_1} |\nabla_\eps\varphi_\eps|^{p_W} \dd x\le c \quad \text{ for }
p_W=\frac{2}{1+4/q_W}\label{eq:conv4}.
\end{align}
\end{lemma}
\begin{proof}
1. Using Corollary~\ref{cor:FE} and the resulting form of $\calF_\eps$ in the upper line 
in \eqref{eq:formFunctionalPoisson} and the Assumption \ref{assu:WgrowthSO3} we find
\[
\int_{\Omega_1}
\dist^2(\nabla_\eps y_\eps M_\eps^{-1}, \SO(3))\det M_\eps\dd x\le C_W\eps^2,
\]
and for $\eps$ small enough by \ref{assu:Prestrain}
\[
\int_{\Omega_1}
\dist^2(\nabla_\eps y_\eps M_\eps^{-1}, \SO(3))\,\mathrm{d}x\le C\eps^2.
\] 
We argue as in \cite[Subsec.~3.1]{Padilla-Garza2022}:
There exists a measurable rotation field $R(x):\Omega_1\to \SO(3)$ such that
\[
\int_{\Omega_1}|\nabla_\eps y_\eps M_\eps^{-1}-R(x)|^2\,\dd x\le c\eps^2,
\]
meaning that $\|\nabla_\eps y_\eps M_\eps^{-1}-R\|_{L^2}\le c\eps$.
Therefore, $\|\nabla_\eps y_\eps M_\eps^{-1}\|_{L^2}\le
\|R\|_{L^2}+c\eps$ holds. This estimate then ensures $\|\nabla_\eps y_\eps M_\eps^{-1}\|_{L^2}\le c$ for small $\eps$. Since 
$B$ is bounded, Assumption~\ref{assu:Prestrain} implies that, for 
$\eps>0$ sufficiently small, $\|\nabla_\eps y_\eps\|_{L^2}\le c$. Next, we estimate by the triangle inequality
\begin{equation*}
\begin{split}
&\int_{\Omega_1}
\dist^2(\nabla_\eps y_\eps, \SO(3))\,\mathrm{d}x\le
2\!\!\int_{\Omega_1}\!
\dist^2(\nabla_\eps y_\eps, \nabla_\eps y_\eps M_\eps^{-1})\,\mathrm{d}x
+2\!\!\int_{\Omega_1}\!
\dist^2(\nabla_\eps y_\eps M_\eps^{-1}, \SO(3))\,\mathrm{d}x\\
&\le c \|\nabla_\eps y_\eps\|_{L^2}^2\,
\esssup_{x\in\Omega_1}\,\dist^2\big(\mathbb{I}_3, M_\eps(x)^{-1}\big)
+c\eps^2\calF_\eps(y_\eps,\varphi_\eps)\leq C\eps^2.
\end{split}
\end{equation*}

2. We write 
$\det \nabla_\eps y_\eps=\det(\nabla_\eps y_\eps M_\eps^{-1})\det M_\eps$ and observe that by Assumption~\ref{assu:Prestrain}, $\det M_\eps>0$ for $\eps$ sufficiently small.
Moreover, since $\calF_\eps(y_\eps,\varphi_\eps)\leq C$ and by the growth condition in Assumption \ref{assu:Wgrowth}, the first factor is positive almost everywhere in $\Omega_1$. Thus $\det \nabla_\eps y_\eps>0$ a.e.\ in $\Omega_1$
for $0<\eps<\eps_0$.

3. Using $|\nabla_\eps y_\eps|\le c |\nabla_\eps y_\eps M_\eps^{-1}|$,  
$|\det (\nabla_\eps y_\eps)|^{-1}\le c 
|\det \nabla_\eps y_\eps M_\eps^{-1}|^{-1}$, and the growth condition in Assumption \ref{assu:Wgrowth},
 we derive 
\[
\int_{\Omega_1}\Big(|\nabla_\eps y_\eps|^{q_W}
+|\det (\nabla_\eps y_\eps)|^{-\frac{q_W}{2}}\Big)\dd x
\le c(q_W,C_W)
\int_{\Omega_1} \Big(W_\text{el}(x,\nabla_\eps y_\eps M_\eps^{-1})+1\Big)\dd x\le c.
\]

4. Furthermore, working with the form \eqref{eq:formFunctionalPoisson} of $\calF_\eps$ we find
\[\int_{\Omega_1} \bfk \big(F_\eps^{-\top}\nabla_\eps \varphi_\eps\big)\cdot \big(F_\eps^{-\top}\nabla_\eps\varphi_\eps\big)\det F_\eps\dd x\le C.
\]
By Assumption~\ref{eq:AssuElectric2},
 $\bfk\in L^{\infty}(\Omega_1,\R^{3\times 3})$ is positive definite, which leads to 
\[\int_{\Omega_1} (F_\eps^{-\top}\nabla_\eps \varphi_\eps(\det F_\eps)^{1/2})\cdot (F_\eps^{-\top}\nabla_\eps\varphi_\eps(\det F_\eps)^{1/2})\dd x
=\int_{\Omega_1} |F_\eps^{-\top}\nabla_\eps \varphi_\eps|^2 \det F_\eps\dd x\le C.
\] 
Since $|\nabla_\eps \varphi_\eps|\le (\det F_\eps)^{1/2}
|F_\eps^{-\top}\nabla_\eps \varphi_\eps||F_\eps|
(\det F_\eps)^{-1/2}$ it follows  from H\"older's inequality with $p_W=\frac{2}{1+4/q_W}$ that
\[
\|\nabla_\eps \varphi_\eps\|_{L^{p_W}}\le
\Big(\int_{\Omega_1}|F_\eps^{-\top}\nabla_\eps \varphi_\eps|^2 \det F_\eps\,\dd x\Big)^{1/2}
\|F_\eps\|_{L^{q_W}} \|(\det F_\eps) ^{-1}\|_{L^{q_W/2}}^{1/2}\le c.
\]
This completes the proof.
\end{proof}

\begin{proposition}\label{prop:conv}
We assume \ref{assu:W} -- \ref{eq:AssuElectric2}.
Let $(y_\eps,\varphi_\eps)_{\eps>0}\subseteq \calY_1 \times \calV_1$ 
be a sequence fulfilling 
 \eqref{eq:le}.
Then, there exist $y_0\in \calY_0$, 
$\varphi_0\in \calV_0$, $m_0\in L^{2}(\omega)$, and a non-relabeled subsequence $(y_\eps,\varphi_\eps)_{\eps>0}$ such that
\begin{align*}
&y_\eps\to y_0\,\text{ in }L^2(\Omega_1,\R^3),
\quad
\nabla_\eps y_\eps\to (
\nabla' y_0|\nu_{y_0})=:R_{y_0}\,\text{ in }L^2(\Omega_1,\R^{3\times 3}),
\\
&\varphi_\eps\to \varphi_0\text{ in }L^2(\Omega_1), \quad
\nabla_\eps\varphi_\eps \rightharpoonup
(\nabla' \varphi_0,m_0)\text{ in }L^{p_W}(\Omega_1,\R^3),
\end{align*}
where $\nu_{y_0}$ denotes the unit normal vector corresponding to $y_0$ (cf.\ \eqref{eq:mu}) and $p_W:=\frac{2}{1+4/q_W}$.
In particular, the convergences  $
y_\eps\to y_0$ in $\calY$,
$\varphi_\eps\to\varphi_0$ in $\calV$ hold true in the norm topology of the spaces 
\begin{equation}\label{eq:spacesYV}
\calY:= W^{1,2}(\Omega_1,\R^3),\quad 
\calV:=L^2(\Omega_1).
\end{equation}
\end{proposition}

\begin{proof}
1. 
The estimates in \eqref{eq:conv1} for the sequence $(y_\eps)_{\eps>0}\subseteq W^{1,2}(\Omega_1,\R^3)$ imply
\[
\limsup_{\eps\to 0} \frac{1}{\eps^2}
\int_{\Omega_1}
\dist^2(\nabla_\eps y_\eps, \SO(3))\,\mathrm{d}x<\infty,
\]
and therefore allow us to apply the compactness result of \cite[Theorem 4.1]{FJM1}. This  ensures that the sequence $(\nabla_\eps y_\eps)_{\eps>0}$ is precompact in $L^2(\Omega_1,\mathbb{R}^{3\times 3})$ 
and that there exist $(\nabla' y_0,b)\in H^1(\Omega_1;\R^{3\times 3})$
and a non-relabeled subsequence such that 
$ \nabla_\eps y_\eps\to (\nabla'y_0|b)$ in $L^2(\Omega_1,\mathbb{R}^{3\times 3})$, where
$(\nabla' y_0|b)\in \SO(3)$ a.e.\ in $\Omega_1$. 
This latter property implies that $b=\nu_{y_0}$ with $\nu_{y_0}$ defined in \eqref{eq:mu}.

Since $\frac{1}{\eps}\partial_3 y_\eps\to \nu_{y_0}$ 
in $L^2(\Omega_1,\mathbb{R}^3)$, we obtain 
$\partial_3 y_\eps\to 0$ in $L^2(\Omega_1,\mathbb{R}^3)$ and the limit $y_0$ does not depend on $x_3$.
Furthermore, $(\nabla'y_0|\nu_{y_0})$ is independent of $x_3$, and 
$(\nabla'y_0|\nu_{y_0})\in W^{1,2}(\omega,\mathbb{R}^{3\times 3}),$ meaning that $(\nabla'y_0|\nu_{y_0})$ is much more regular than naively expected.

According to the estimate in \eqref{eq:conv2}, we can conclude $\|\nabla_\eps y_\eps\|_{L^{q_W}}\le c$  and therefore $\|\nabla y_\eps\|_{L^2}\le c$.
As $\int_{\Omega_1}y_\eps\dd x=0$, we have by the Poincar\'e--Wirtinger inequality
$\|y_\eps\|_{W^{1,2}}\le c$. Thus, 
$(y_\eps)_{\eps>0}$ is precompact in $L^2(\Omega_1;\R^3)$ and for a 
non-relabeled subsequence we obtain 
$y_\eps\to y_0$ in $L^2(\Omega_1,\R^3)$ with the same limit $y_0$ as above.

2.
From \eqref{eq:conv4} and $\varphi_\eps\in \calV_1$, we obtain
 $\|\varphi_\eps\|_{W^{1,p_W}(\Omega_1)}\le c$.
Therefore, we find 
$\varphi_0\in  W^{1,p_W}(\Omega_1)$ such that for a subsequence
$\varphi_\eps\rightharpoonup \varphi_0$ in $W^{1,p_W}(\Omega_1)$,
 and for $q_W>6$ (cf.\ \ref{assu:Wgrowth}, implying $p_W>6/5$) 
 $\varphi_\eps\to \varphi_0$ in $L^2(\Omega_1)$ by compact embedding.
Moreover, $\nabla_\eps\varphi_\eps \rightharpoonup (\nabla'\varphi_0,m_0)$ in 
$L^{p_W}(\Omega_1,\R^3)$  can be assumed with some $m_0\in L^{p_W}(\Omega_1)$.
From $\tfrac{1}{\eps}\partial_3\varphi_\eps\rightharpoonup m_0$ in $L^{p_W}(\Omega_1)$, we find that $\partial_3\varphi_\eps\to 0$ in $L^{p_W}(\Omega_1)$. Therefore, the limit $\varphi_0$ does not depend on $x_3$, and we can identify it with a function 
$\varphi_0\in W^{1,p_W}(\omega)$. 

Since we only have that $6/5<p_W< 2$, 
our next aim is to show that
$(\nabla' \varphi_0,m_0)\in L^2(\Omega_1,\R^3)$. We follow the ideas 
in \cite[pp.\ 1483--1484]{ag-bar23} and show the following auxiliary result:

Let $Z_\eps \rightharpoonup Z$ in $L^{p_W}(\Omega_1,\R^3)$,
$F_\eps\to R_{y_0}$ in $L^2(\Omega_1,\R^{3\times 3})$ such that 
we have $\det F_\eps>0$ and $R_{y_0}\in \SO(3)$ a.e.\ in $\Omega_1$. Then, it holds that
\begin{equation}\label{eq:hilfa}
\int_{\Omega_1}|Z|^2\,\dd x
\le \liminf_{\eps\to 0}  
\int_{\Omega_1}|F_\eps^{-\top}Z_\eps|^2 \det F_\eps\,\dd x.
\end{equation}

For the proof of \eqref{eq:hilfa}, we use the measurable and bounded map
\[
\Psi:\R^{3\times 3}\to \R^{3\times 3},\quad
\Psi(F):=\begin{cases}
F^{-\top}\sqrt{\det F}& \text {if }\mathrm{dist}(F,\SO(3))\le \frac{1}{2},\\
0 & \text{else}
\end{cases}
\] 
with $C_\Psi:=\sup_{F\in\R^{3\times 3}}|\Psi(F)|<\infty$.
For a subsequence, we have $F_\eps \to R_{y_0}$ a.e.\ in $\Omega_1$. 
The continuity of $\Psi$ in an open neighborhood of $\SO(3)$
and the property $R_{y_0}\in \SO(3)$ a.e.\ in $\Omega_1$ then ensure 
$\Psi(F_\eps)\to \Psi(R)=R_{y_0}$ a.e.\ in $\Omega_1$.
Since $\Psi(F_\eps)\le C_\Psi$, we conclude from the weak convergence of $Z_\eps$ in $L^{p_W}(\Omega_1,\R^3)$ that also
$\Psi(F_\eps)Z_\eps\rightharpoonup R_{y_0}Z$ in $L^{p_W}(\Omega_1,\R^3)$ holds. Therefore, the weak lower semicontinuity of convex  functionals yields
\[
\int_{\Omega_1}|R_{y_0}Z|^2\,\dd x\le \liminf_{\eps \to 0}
\int_{\Omega_1}|\Psi(F_\eps)Z_\eps|^2\,\dd x.
\]
Using that pointwise $|RZ|=|Z|$ and 
$|\Psi(F_\eps)Z_\eps|^2\le |F_\eps^{-\top} Z_\eps|^2|\det F_\eps|$
(by  the definition of $\Psi$), we obtain the desired estimate \eqref{eq:hilfa}.

Applying now \eqref{eq:hilfa} for 
$Z=(\nabla'\varphi_0,m_0)$, 
$Z_\eps=\nabla_\eps\varphi_\eps$, and $F_\eps=\nabla_\eps y_\eps$ together with the estimate \eqref{eq:conv3} gives $(\nabla' \varphi_0,m_0)\in L^2(\Omega_1,\R^3)$.
From $\varphi_\eps\in L^2_\text{av}(\Omega_1)$ and 
$\varphi_\eps\to \varphi_0$ in $L^2(\Omega_1)$, we finally conclude $\varphi_0\in \calV_0$.
\end{proof}

\section{Main result: dimension reduction for the electro-elastic problem}
\label{sec:dimred}

This section contains the main result of the paper, namely, the convergence of the bivariate functionals $\calF_\eps$ defined in \eqref{eq:totalEnergyScaled}. First, we 
introduce our abstract notion of convergence building upon \cite{AttouchWets83}. The dimension reduction is carried out in Subsection~\ref{ss:mainresult}.

\subsection{Abstract convergence result for saddle point problems}
We first introduce the abstract framework that will be used in the sequel, building on the notion of epi/hypo-convergence developed in \cite{AW1981ACNO,AttouchWets83}.
Let $\calY, \calY_0,\calY_1$ and $\calV, \calV_0,\calV_1$ be metric spaces such that $\calY_i\subseteq\calY$ and $\calV_i\subseteq\calV$ for $i=0,1$. 
Consider a family of bivariate functionals 
$\calF_n:\calY_1\times\calV_1\to \R\cup\{\infty\}$ and a  limit functional $\calF_\infty:\calY_0\times\calV_0\to\R$.

A saddle point of $\calF_n$ is a pair $(y_n^*,\varphi_n^*)\in\calY_1\times\calV_1$ such that
\begin{equation}\label{eq:abstractSaddleProblemFamily}
\calF_n(y_n^*,\wh \varphi)
\leq \calF_n(y_n^*,\varphi_n^*)\leq \calF_n(\wh y,\varphi_n^*)\quad\forall (\wh y,\wh \varphi)\in \calY_1\times\calV_1.
\end{equation}
Similarly, a saddle point of $\calF_\infty$ is a pair $(y_0^*,\varphi_0^*)\in\calY_0\times\calV_0$ such that
\begin{equation}\label{eq:abstractSaddleProblemFamilyLimit}
\calF_\infty(y_0^*,\wh \varphi)
\leq \calF_\infty(y_0^*,\varphi_0^*)\leq \calF_\infty(\wh y,\varphi_0^*)\quad\forall (\wh y,\wh \varphi)\in \calY_0\times\calV_0.
\end{equation}

\begin{proposition}
\label{prop:abstractResult}
    Consider a family of bivariate functionals 
$\calF_n:\calY_1\times\calV_1\to \R\cup\{\infty\}$ with saddle points $(y_n^*,\varphi_n^*)\in\calY_1\times\calV_1$
with $\sup_n\calF_n(y_n^*,\varphi_n^*)<\infty$ and
such that $(y_n^*,\varphi_n^*)\to(y_0^*,\varphi_0^*)\in \calY_0\times \calV_0$ with convergence in $\calY\times\calV$. Assume that 
\begin{align}
\label{eq:EHliminf}
\forall\,\wh \varphi_0\in\calV_0 \ \exists\,(\wh \varphi_n)\subseteq\calV_1 &:\ 
\wh \varphi_n\to \wh \varphi_0\text{ in }\calV\text{ and }\liminf_{n\to\infty} \calF_n(y_n^*,\wh \varphi_n)\geq \calF_\infty(y_0^*,\wh \varphi_0),\\
\label{eq:EHlimsup}
\forall\,\wh y_0\in\calY_0 \ \exists\,(\wh y_n)\subseteq\calY_1&:\ 
\wh y_n\to \wh y_0\text{ in }\calY\text{ and }\limsup_{n\to\infty}\calF_n(\wh y_n,\varphi_n^*)\leq \calF_\infty(\wh y_0,\varphi_0^*).
\end{align}
Then $(y_0^*,\varphi_0^*)$ is a saddle point for $\calF_\infty$ and
$\lim_{n\to\infty}\calF_n(y_n^*,v_n^*)=\calF_\infty(y_0^*,v_0^*)$.
\end{proposition}

\begin{proof}
Let us denote $a_*:=\liminf_{n\to\infty}\calF_n(y_n^*,\varphi_n^*)\leq \limsup_{n\to\infty}\calF_n(y_n^*,\varphi_n^*)=:b_*$. Let $\wh \varphi_0\in\calV_0$ be arbitrary. We choose a sequence $(\wh \varphi_n)\subseteq \calV_1$ such that 
$\wh\varphi_n \to \wh\varphi_0$ in $\calV$ as in \eqref{eq:EHliminf}.   Therefore, by \eqref{eq:abstractSaddleProblemFamily} we obtain
\[\calF_\infty(y^*_0,\wh \varphi_0)\leq \liminf_{n\to\infty}\calF_n(y_n^*,\wh\varphi_n)\leq\liminf_{n\to\infty}\calF_n(y_n^*,\varphi_n^*)=a_*.\]
Similarly, for arbitrary $\wh y_0\in\calY_0$ we choose a sequence $(\wh y_n)\subseteq \calY_1$ with $\wh y_n\to \wh y_0$ in $\calY$ satisfying \eqref{eq:EHlimsup} to get $b_* \leq \calF_\infty(\wh y_0,\varphi^*_0)$.
This holds for any pairs $(\wh y_0,\wh \varphi_0)\in\calY_0\times\calV_0$. Thus, it results
\begin{eqnarray}\label{eq:limsaddle}
\calF_\infty(y^*_0,\wh\varphi_0)\leq a_*\leq b_*\leq \calF_\infty(\wh y_0,\varphi^*_0)\quad\text{for all }(\wh y_0,\wh \varphi_0)\in\calY_0\times\calV_0.
\end{eqnarray}
Finally, choosing $\wh \varphi_0 =\varphi^*_0$ in \eqref{eq:limsaddle} gives
$\calF_\infty(y^*_0,\varphi^*_0)\leq\calF_\infty(\wh y_0,\varphi^*_0)$ for all $\wh y_0\in\calY_0$, and analogously for $\wh y_0=y^*_0$
we obtain $\calF_\infty(y^*_0,\varphi^*_0) \geq\calF_\infty(y^*_0,\wh\varphi_0)$ for all $\wh\varphi_0\in\calV_0$. Thus, $(y_0^*,\varphi_0^*)$ is a saddle point for $\calF_\infty$. \\

It remains to show the convergence $\calF_n(y_n,\varphi_n)\to\calF_\infty(y^*_0,\varphi^*_0)$. Indeed, by \eqref{eq:EHliminf} there exists a sequence $(\wh \varphi_n)\subseteq\calV_1$ such that $\wh \varphi_n\to\varphi^*_0$ in $\calV$ and 
\[
 \calF_\infty(y^*_0,\varphi^*_0)\leq\liminf_{n\to\infty}\calF_n(y_n^*,\wh\varphi_n)\leq\liminf_{n\to\infty}\calF_n(y_n^*,\varphi_n^*)=a_*,
\]
where we used that $(y_n^*,\varphi_n^*)$ is a saddle point.
Analogously, we get from \eqref{eq:EHlimsup} that $b_*\leq \calF_\infty(y^*_0,\varphi^*_0)$.
Hence, $a_*=b_*=\calF_\infty(y^*_0,\varphi^*_0)$.
\end{proof}

\begin{remark}
The theory of epi/hypo-convergence of saddle-point problems developed in \cite{AttouchWets83} (see also \cite{AW1981ACNO}) is formulated only for one space  $\calY\times\calV$, 
instead of the triple $\calY_0\times\calV_0,\calY\times\calV,\calY_1\times\calV_1$ as in our more specific case. 
Moreover, the upper and lower estimates in \eqref{eq:EHliminf} and \eqref{eq:EHlimsup} are formulated for general sequences instead of saddle points.
\end{remark}

\subsection{Dimension reduction}
\label{ss:mainresult}

We aim to apply the abstract result of Proposition~\ref{prop:abstractResult} to the 
electro-elastic saddle-point problem for the functional \eqref{eq:totalEnergyScaled} using the spaces 
$\calY_0,\,\calV_0$ introduced in \eqref{eq:spacesY0V0},
$\calY_1,\,\calV_1$ defined in \eqref{eq:spacesY1V1}, 
and $\calY,\,\calV$ specified in \eqref{eq:spacesYV}. 
In particular, we show that for $\eps\to 0$ the $\liminf$ and $\limsup$ conditions
\eqref{eq:EHliminf} and \eqref{eq:EHlimsup} are satisfied 
(meaning that for every sequence $\eps_n\to 0$, as  $n\to \infty$, 
the conditions hold for $\calF_n = \calF_{\eps_n}$).

We assume in the following that 
 \begin{equation}\label{eq:le}
 \begin{gathered}
 (y_\eps,\varphi_\eps)_{\eps>0}\subseteq \calY_1 \times \calV_1 
 \text{ is a sequence such that }\\
 \calF_\eps(y_\eps,\wh\varphi)\le 
 \calF_\eps(y_\eps,\varphi_\eps)~\forall 
 \wh\varphi\in   
 \calV_1, \text{ and }\exists\, C>0:\quad
 \calF_\eps(y_\eps,\varphi_\eps)\leq C.
 \end{gathered}
 \end{equation}
In particular, Lemma~\ref{lem:apri} and Proposition \ref{prop:conv} apply, and we find (non-relabeled) subsequences and limits $(y_0,\varphi_0)\in\calY_0\times \calV_0$ such that $y_\eps\to y_0$ in $\calY$, $\varphi_\eps\to \varphi_0$ in $\calV$. Recall the definitions of the limiting functionals $\calM_0,\calE_0 $ and $\calF_0$ in Subsection~\ref{ss:limitModel}. We proceed 
in the following steps: First, we show that
\begin{align*}
&\text{(I)}~
\liminf_{\eps\to 0} \calM_\eps(y_\eps)\geq \calM_0(y_0),
\quad
\\
\nonumber
&\text{(II)}~
\forall\, \wh \varphi_0\in \mathcal{V}_0\, \exists (\wh \varphi_\eps)_{\eps>0}\subseteq \calV_1:\quad
\wh \varphi_\eps\to\wh \varphi_0\text{ in $\calV$\ and\ }\limsup_{\eps\to 0} \calE_\eps(y_\eps,\wh\varphi_\eps) \leq \calE_0(y_0,\wh\varphi_0).
\end{align*}
If (I) and (II) are satisfied, the $\liminf$ condition in \eqref{eq:EHliminf} also holds since
\[
\liminf_{\eps\to 0}\calF_\eps(y_\eps,\wh \varphi_\eps)
\geq \liminf_{\eps\to 0}\calM_\eps(y_\eps) - \limsup_{\eps\to 0}\calE_\eps(y_\eps,\wh \varphi_\eps)
\geq \calM_0(y_0) -\calE_0(y_0,\wh\varphi_0)=\calF_0(y_0,\wh \varphi_0).
\]
Next, we show that the $\limsup$ condition in \eqref{eq:EHlimsup} also holds, i.e.,
\begin{align*}
\nonumber
&\text{(III)}~
\forall\, 
\wh y_0\in \mathcal{Y}_0\,\exists \,(\wh y_\eps)_{\eps>0}\subseteq \calY_1:\quad\wh y_\eps\to \wh y_0\text{ and }
 \limsup_{\eps\to 0} \calF_\eps(\wh y_\eps,\varphi_\eps) \leq \calF_0(\wh y_0,\varphi_0).
\end{align*}
Finally, while not necessary for Proposition~\ref{prop:abstractResult}, we also show that
\begin{align}
&\text{(IV)}~\liminf_{\eps\to 0} \calE_\eps(y_\eps,\varphi_\eps)\geq \calE_0(y_0,\varphi_0).
\end{align}

\begin{remark}
\label{rem:reformElectrostatic}
The conditions (II) and (IV) imply the $\Gamma$-convergence of the functionals
$\wt \calE_\eps(\cdot) := \calE_\eps(y_\eps,\cdot)$ for the fixed sequence of deformations $y_\eps$.
We will show that the sequence $\wh y_\eps$ in (III) is also a recovery sequence for $\calM_\eps$ such that
together with (I) the $\Gamma$-convergence of the mechanical energy $\calM_\eps$ also follows.
\end{remark}

We introduce the 
following additional assumption 
as in \cite{ag-ago19}
\begin{enumerate}[label=(A6)]
\item \label{assu:B} 
$W_\mathrm{el}(x',x_3)=W_\mathrm{el}(x')\quad\text{and $B(x)$ is such that}\quad
\mathrm{curl} \Big(\mathrm{curl} \displaystyle\int_{-1/2}^{1/2}B_{2\times 2}(x',t)\dd t\Big)=0$.
\end{enumerate}
The operator 'curl' inside
the parenthesis acts on a $2\times 2$ matrix by taking the 'curl' of each row, giving as a result a
two-dimensional vector.
Note that the necessity of this assumption was removed
in \cite{Padilla-Garza2022}. For the sake of clarity and comprehensibility in our proofs, 
we have kept the assumption. See also Section~\ref{sec:conclusion}.

The Assumption~\ref{assu:B} ensures by \cite[Theorem 2.8]{ag-ago19}, resp. \cite[Theorem 3.2]{ag-cia05} that there exists a vector potential
$g\in W^{1,2}(\omega,\R^2)$ such that
\begin{equation}\label{vektorpotential}
\int_{-1/2}^{1/2} 
B_{2\times 2}(x',t)\dd t=\nabla_{\mathrm{sym}}g
:=\frac{1}{2}(\nabla' g+\nabla' g^\top),
\end{equation}
where $g$ is unique up to rigid displacements. 

Now we are ready to formulate the main result 
of our paper:

\begin{theorem}\label{thm:main}
We assume \ref{assu:W} -- \ref{eq:AssuElectric2}.
Let $(y_\eps,\varphi_\eps)_{\eps>0}\subseteq \calY_1 \times \calV_1$ 
be a sequence as in \eqref{eq:le}
Then there exist $y_0\in \calY_0$, 
$\varphi_0\in \calV_0$, $m_0\in L^{2}(\omega)$, and a non-relabeled subsequence $(y_\eps,\varphi_\eps)_{\eps>0}$ such that
\begin{align*}
&y_\eps\to y_0\text{ in }L^2(\Omega_1,\R^3),
\quad
\nabla_\eps y_\eps\to (
\nabla' y_0|\nu_{y_0})=:R_{y_0}\text{ in }L^2(\Omega_1,\R^{3\times 3}),
\\
&\varphi_\eps\to \varphi_0\text{ in }L^2(\Omega_1), \quad
\nabla_\eps\varphi_\eps \rightharpoonup
(\nabla' \varphi_0,m_0)\text{ in }L^{p_W}(\Omega_1,\R^3),
\end{align*}
where $\nu_{y_0}$ denotes the unit normal vector corresponding to $y_0$ and $p_W:=\frac{2}{1+4/q_W}$.
Moreover, the conditions {\rm (I)}, {\rm (II)}, and {\rm (IV)} are satisfied. 
Supposing additionally \ref{assu:B}, then also Condition {\rm (III)} is fulfilled.
\end{theorem}

The proof of Theorem~\ref{thm:main} results from Lemma~\ref{lem:apri} and Proposition \ref{prop:conv} in Subsection~\ref{ss:boundedness} and Lemmas~\ref{lem:I},
\ref{lem:II}, \ref{lem:III}, and \ref{lem:IV} in Subsection~\ref{subsec26}.

\begin{corollary}We assume \ref{assu:W} -- \ref{assu:B}.
Then the saddle-point problems for the functionals $\calF_\eps$ converge to the saddle-point problem associated with the two-dimensional limit functional $\calF_0$ in the sense of Proposition~\ref{prop:abstractResult}
and limits of saddle points $(y_\eps^*,\varphi_\eps^*)$ for $\calF_\eps$ are saddle points for $\calF_0$.
\end{corollary}

We highlight that Theorem~\ref{thm:main} is proven for sequences $(y_\eps,\varphi_\eps)$
satisfying~\eqref{eq:le}, while the above corollary considers  sequences $(y_\eps^*,\varphi_\eps^*)$ that are additionally saddle points.

\subsection{Proof of asymptotic upper and lower bounds}\label{subsec26}

We start with property  (I):

\begin{lemma}\label{lem:I}{\rm{[Property (I)].}}
We assume \ref{assu:W} -- \ref{eq:AssuElectric2}.
Let $(y_\eps,\varphi_\eps)_{\eps>0}\subseteq \calY_1 \times \calV_1$ 
be a sequence fulfilling \eqref{eq:le}. Moreover, let $y_0\in \calY_0$  be a 
limit  with respect to the convergence stated in Proposition~\ref{prop:conv}. Then, it holds 
\[
\liminf_{\eps\to 0} \calM_\eps(y_\eps)\geq \calM_0(y_0).
\]
\end{lemma}

\begin{proof}
The non-negativity of the hyper stress term 
(see \ref{assu:hyper}) ensures
\begin{equation*}
\begin{split}
\calM_\eps(y_\eps)
&\ge 
 \frac{1}{\eps^2}\int_{\Omega_1} 
    W_\text{el}( \cdot,\nabla_\eps  y_\eps   M_\eps^{-1})\det( M_\eps)\dd x\\
&= 
\frac{1}{\eps^2}\int_{\Omega_1}
W_\text{el}( \cdot,\nabla_\eps  y_\eps   M_\eps^{-1})(\det( M_\eps)-1)\dd x
+\frac{1}{\eps^2}\int_{\Omega_1} 
    W_\text{el}( \cdot,\nabla_\eps  y_\eps   M_\eps^{-1})\dd x.
\end{split}
\end{equation*}
Since  $\det(M_\eps)=1+\eps\mathrm{tr}B+O(\eps^2)$ (see Remark~\ref{rem:assu}), 
we find $\|\det( M_\eps){-}1\|_{L^\infty}\to 0$ for $\eps\to 0$, which together with the uniform boundedness of $\frac{1}{\eps^2}\int_{\Omega_1} 
    W_\text{el}( \cdot,\nabla_\eps  y_\eps   M_\eps^{-1})\dd x$
    (for $\eps$ small enough) implies that the first term on the second line tends to zero
    for $\eps\to 0$.
As in the proof of \cite[Theorem 1, Eq.~(23), esp. Eq.~(61)]{Padilla-Garza2022}, we obtain for the second term  $\liminf_{\eps\to 0}\frac{1}{\eps^2}\int_{\Omega_1} 
    W_\text{el}( \cdot,\nabla_\eps  y_\eps   M_\eps^{-1})\dd x
\ge \calM_0(y_0)$, which finishes the proof.
\end{proof}

We prove (II):
\begin{lemma}\label{lem:IV}{\rm{[Property (II)].}}
We assume \ref{assu:W} -- \ref{eq:AssuElectric2}.
Let $(y_\eps,\varphi_\eps)_{\eps>0}\subseteq \calY_1 \times \calV_1$ 
be a sequence fulfilling \eqref{eq:le}.
Moreover, let $y_0\in \calY_0$ be the limit 
with respect to the convergence stated in Proposition~\ref{prop:conv}. 
Then for all $\wh \varphi\in \calV_0$ there exists a sequence 
$(\wh \varphi_\eps)_{\eps>0}\subseteq \calV_1$ such that $\wh \varphi_\eps\to \wh\varphi$ in  $L^2(\Omega_1)$ and
\[
\limsup_{\eps\to 0} \calE_\eps(y_\eps,\wh \varphi_\eps)
= \calE_0(y_0,\wh \varphi).
\]
\end{lemma}

\begin{proof}
Let $\widehat \varphi\in \calV_0$ be arbitrary given and let $m=m_{y_0,\widehat \varphi}\in L^2(\omega)$ be as in the definition of $\calE_0(y_0,\wh \varphi)$, see  \eqref{eq:defMYPHI}. 
Since $C^\infty(\overline{\omega})$ is dense in  $\calV_0$ and $L^2(\omega)$, we find for any $l\in\mathbb{N}$ some $C^\infty(\overline{\omega})$ functions $\wh \varphi_l$, $m_l$ such that
\[
\|\wh \varphi-\wh \varphi_l\|_{W^{1,2}(\omega)}\le \frac{1}{4l},\quad 
\|\wh \varphi-\wh \varphi_l\|_{L^1(\omega)}\le \frac{1}{4l\gamma\|n_\mathrm{ch}\|_{L^\infty(\Omega_1)}},\quad
\|m-m_l\|_{L^2(\omega)}\le \frac{1}{4l},
\quad\text{ and }
\]
\[
\frac{\beta}{2}\int_{\Omega_1}
\Big|R_{y_0}^\top\bfk R_{y_0} (\nabla' \widehat \varphi_l,m_l)\cdot (\nabla' \widehat \varphi_l,m_l)
-
R_{y_0}^\top\bfk R_{y_0} (\nabla' \widehat \varphi,m)\cdot (\nabla' \widehat \varphi,m)\Big|\,\dd 
x\le \frac{1}{4l}.
\]
We define
\[
\widehat \varphi_{l\eps} (x',x_3):=\widehat \varphi_l(x')+\eps m_l(x') x_3,\quad\text{ thus \quad}
\widehat \varphi_{l\eps}\in C^\infty(\overline{\Omega_1}),
W^{1,\infty}(\Omega_1). 
\]
From this definition, it follows directly that
\[
\widehat \varphi_{l\eps}\to \widehat \varphi_l,\,\,
\nabla_\eps \widehat \varphi_{l\eps}
=(\nabla'\widehat \varphi_l+\eps x_3\nabla'm_l , m_l)\to  
(\nabla'\widehat \varphi_l, m_l)\text{ in } 
L^\infty(\Omega_1)\text{ as }\eps\to 0 .
\]

For the given sequence $(y_\eps)_{\eps>0}$, let $F_\eps:=\nabla_\eps y_\eps$, 
we apply 
Lemma~\ref{lem:qw6} below to verify that 
\begin{equation*}
\begin{split}
\lim_{\eps\to 0}\frac{\beta}{2}\int_{\Omega_1} \!\!
F_\eps^{-1}\bfk F_\eps^{-\top}(\nabla_\eps \widehat \varphi_{l\eps})\cdot (\nabla_\eps \widehat \varphi_{l\eps} )\det F_\eps\,\dd x
=
\frac{\beta}{2}\int_{\Omega_1}\!\! R_{y_0}^\top\bfk R_{y_0} (\nabla' \widehat \varphi_l,m_l)\cdot (\nabla' \widehat \varphi_l,m_l)\,\dd x
.
\end{split}
\end{equation*}	
For $l\in\mathbb{N}$, we choose a monotonously decaying sequence $\eps_l=\eps(l)$, 
	$\eps(l)<\eps(l{-}1)$, such that $\eps_l\downarrow 0$ and for all $\eps\le \eps_l$
\begin{equation*}
\begin{split}
&\Big|\int_{\Omega_1}\!\! R_{y_0}^\top\bfk R_{y_0} (\nabla' \widehat \varphi_l,m_l)\cdot (\nabla' \widehat \varphi_l,m_l)\,\dd x
-
\int_{\Omega_1} \!\!
F_{\eps}^{-1}\bfk F_{\eps}^{-\top}(\nabla_{\eps} \widehat \varphi_{l\eps})\cdot (\nabla_{\eps} \widehat \varphi_{l\eps} )\det F_{\eps}\,\dd x\Big|\le \frac{1}{2l\beta},\\
&\text{and }
\|\wh \varphi_{l\eps_l}-\wh \varphi_l\|_{L^1(\Omega_1)}\le \frac{1}{4l\gamma\|n_\mathrm{ch}\|_{L^\infty(\Omega_1)}}.
\end{split}
\end{equation*}
For $\eps\in(0,\eps_0)$, we further define the final sequence
\[
\wh \varphi_\eps:=\wh \varphi_{l\eps}\quad\text{ if }\eps\in(\eps_{l+1},\eps_l]
\text{ for this }l\in\mathbb{N}.
\]
By the previous estimates, we obtain 
for $\eps\in(\eps_{l+1},\eps_l]$
\begin{equation*}
\begin{split}
&|\calE_0(y_0,\widehat \varphi)-\calE_{\eps}(y_{\eps},\widehat \varphi_{\eps})|\\
&\le\gamma\Big|\int_{\Omega_1}n_\mathrm{ch}(\wh\varphi-\wh\varphi_{\eps})\dd x\Big|\\
&+\frac{\beta}{2}\Big|\int_{\Omega_1}\big\{
F_{\eps}^{-1}\bfk F_{\eps}^{-\top}(\nabla_{\eps} \widehat \varphi_{\eps})\cdot (\nabla_{\eps} \widehat \varphi_{\eps} )\det F_{\eps}
-R_{y_0}^\top\bfk R_{y_0} (\nabla' \widehat \varphi_l,m_l)\cdot (\nabla' \widehat \varphi_l,m_l)\big\}\dd 
x\,\Big|\\
&+\frac{\beta}{2}\Big|\int_{\Omega_1}\big\{
R_{y_0}^\top\bfk R_{y_0} (\nabla' \widehat \varphi_l,m_l)\cdot (\nabla' \widehat \varphi_l,m_l)
-
R_{y_0}^\top\bfk R_{y_0} (\nabla' \widehat \varphi,m)\cdot (\nabla' \widehat \varphi,m)\big\}\dd 
x\,\Big|\le \frac{1}{l},
\end{split}
\end{equation*}		
$l\in \mathbb{N}$.
Thus, letting $l\to \infty,$ meaning 
$(\eps_{l+1},\eps_l]\ni\eps\downarrow 0,$ we find
$\calE_{\eps}(y_{\eps},\widehat \varphi_{\eps})\to 
\calE_0(y_0,\widehat \varphi)$ as $\eps\to 0$.
This concludes the proof.
\end{proof}

We now turn to  property (III):\\

\begin{lemma}\label{lem:II}{\rm{[Property (III)].}}
In addition to the Assumptions \ref{assu:W} -- \ref{eq:AssuElectric2}, we suppose that Assumption \ref{assu:B} holds,
(comp.\ with \cite[Theorem 2.6]{ag-ago19}).
 Let $(y_\eps,\varphi_\eps)_{\eps>0}\subseteq \calY_1 \times \calV_1$ 
be a sequence fulfilling \eqref{eq:le}.  
 Moreover, let $\varphi_0\in \calV_0$ be the limit 
with respect to the convergence stated in Proposition~\ref{prop:conv}. 
Then, for all $\wh y\in \calY_0$ there exists a sequence 
$(\wh y_\eps)_{\eps>0}\subseteq \calY_1$ such that
$\wh y_\eps\to \wh y$ in $W^{1,2}(\Omega_1,\R^3)$
and
\[
\limsup_{\eps\to 0} \calF_\eps(\wh y_\eps, \varphi_\eps)
\le \calF_0(\wh y_0, \varphi_0).
\]
 \end{lemma}

\begin{proof}
1. Let $\wh y\in \calY_0$ be arbitrarily given. Since functions in $\calY_0\cap C^\infty(\overline\omega,\R^3)$  are dense in $\calY_0$ (see \cite{Hornung11} and \cite{Pakzad04}), we find for $l\in \mathbb{N}$ a sequence  $\wh y_l\in C^\infty(\overline \omega,\R^3)$ such that 
\[
\|\wh y-\wh y_l\|_{W^{2,2}}\le \frac{1}{4l}, \quad
\|(\nabla' \wh y_l|\nu_{\wh y_l})-(\nabla' \wh y|\nu_{\wh y}) \|_{L^2}\le \frac{1}{4l}, \quad
|\calM_0(\wh y_l)-\calM_0(\wh y)|\le \frac{1}{4l}.
\]
This construction, together with $Q_2$ being a quadratic form on $\R^{2\times 2}$, 
$\wh y_l\to \wh y$ in $\calY_0$, and the definition of $\overline{Q}_2$ yields
\[
\calM_0(\wh y)-\calM_0(\wh y_l)
=\frac{1}{2}\int_{\Omega_1}\Big\{\overline{Q}_2(x',\nabla \wh y^\top\nabla \nu_{\wh y}
)
-\overline{Q}_2(x',\nabla \wh y_l^\top\nabla \nu_{\wh y_l}) \Big\}\dd x'
\to 0\quad \text{ for }l\to \infty.
\]

2. For each of the $\wh y_l\in \calY_0$, Theorem~\ref{thm:dicht}  ensures the existence of a sequence of functions in $W^{2,q_H}(\Omega_1,\R^3)$, denoted $\wh y_{l\eps}$, such that $\nabla_\eps \wh y_{l\eps}\to (\nabla' \wh y_l|\nu_{\wh y_l})$ in $L^2(\Omega_1,\R^{3\times 3})$ and 
\[
\lim_{\eps\to 0} \frac{1}{\eps^2}\int_{\Omega_1} 
    W_\text{el}( x',\nabla_\eps  \wh y_{l\eps} ( x)  M_\eps( x)^{-1})\dd x
    =\frac{1}{2}\int_{\Omega_1}\overline{Q}_2(x',\nabla \wh y_l^\top\nabla \nu_{\wh y_l})\,\dd x'.
\]
Assumption \ref{assu:Prestrain} and Remark~\ref{rem:assu} therefore guarantee
\[
\lim_{\eps\to 0} \frac{1}{\eps^2}\int_{\Omega_1} 
    W_\text{el}( \cdot,\nabla_\eps  \wh y_{l\eps}   M_\eps^{-1})
    \det M_\eps \dd x	=
\frac{1}{2}\int_{\Omega_1}\overline{Q}_2(x',\nabla \wh y_l^\top\nabla \nu_{\wh y_l})\,\dd x'.    	
\]
The functions $\wh y_{l\eps}$ in the proof of Theorem~\ref{thm:dicht} have the form (see \eqref{ansatz}) 

\begin{equation}\label{eq:approx}
\wh y_{l\eps}(x',x_3)=\wh y_l(x')
+\eps \big[ x_3 \nu_{\wh y_l}(x')+\nabla' \wh y_l(x')g_\eps(x')\big]
+\eps^2D_{l\eps}(x',x_3),
\end{equation}
where
\[
D_{l\eps}(x',x_3):=\int_0^{x_3}d_{l\eps}(x',t)\dd t
\]
with an $\eps$ dependent regularization $d_{l\eps}$ of a function
$d_l\in L^2(\Omega_1,\R^3)$ 
and  $g_{\eps}$ is an $\eps$  dependent regularization of the vector potential $g$ defined by \eqref{vektorpotential}, (see also Lemma~\ref{mollification}).

3. To use the deformation $\wh y_{l\eps}$ in the energy functional $\calF_\eps$, we  additionally need  $\wh y_{l\eps}\in W^{2,q_H}(\Omega_1,\R^3)$. However, this property is guaranteed by Theorem~\ref{thm:dicht} below. Moreover, we aim to show that 
\begin{equation}\label{eq:regula}
\lim_{\eps\to 0} \frac{1}{\eps^2}\int_{\Omega_1} 
\eps^{\alpha_H}H_*(\nabla_\eps^2\wh y_{l\eps})\dd x= 0.
\end{equation}
The proof of Theorem~\ref{thm:dicht}
shows that
$d_{l\eps}\in W^{2,q_H}(\Omega_1,\R^3)$, 
$g_\eps\in W^{2,q_H}(\omega,\R^2)$ for all $\eps$, as well as
$\|d_{l\eps}\|_{W^{2,q_H}(\Omega_1)}^{q_H}\le c/\eps$ and
$\|g_{\eps}\|_{W^{2,q_H}(\omega)}^{q_H}\le c/\eps$.
Moreover, we obtain $d_{l\eps} \to 
 d_l$  in $L^2(\Omega_1,\R^3)$,
 $g_{\eps} \to 
 g$  in $W^{1,2}(\omega,\R^2)$
  for $\eps\to 0$.
We use Assumption~\ref{assu:hyper},  \eqref{eq:approx},
and \eqref{convdicht} to estimate 
\begin{equation*}
\begin{split}
&\frac{1}{\eps^2}\int_{\Omega_1} 
\eps^{\alpha_H}H_*(\nabla_\eps^2\wh y_{l\eps})\dd x
\le 
\int_{\Omega_1}\eps^{\alpha_H-2}K_H
(1+|\nabla_\eps^2\wh y_{l\eps}|^{q_H})\dd x 
\le 
c\eps^{\alpha_H-2-2q_H}(1+\|\wh y_{l\eps}\|_{W^{2,q_H}}^{q_H})\\
&\le
c\eps^{\alpha_H-2-2q_H}\big(1+\|\wh y_l\|_{C^\infty}^{q_H}
+\eps^{q_H}\big(\|\wh y_l\|_{C^\infty}^{q_H}
+\|\wh y_l\|_{C^\infty}^{q_H}\|g_\eps\|_{W^{2,q_H}}^{q_H}\big)
+\eps^{2q_H}\| d_{l\eps}\|_{W^{2,q_H}}^{q_H}\big)\\
&\le \eps^{\alpha_H-2-2q_H}c(l). 
\end{split}
\end{equation*}
(Note that $\| D_{l\eps}\|_{W^{2,q_H}}^{q_H})\le c \| d_{l\eps}\|_{W^{2,q_H}}^{q_H}$.)
Using  that $\alpha_H>2{+}2q_H$,
we find the desired convergence \eqref{eq:regula}.

4. We choose now a monotonically  decreasing subsequence $\eps_l=\eps(l)$ such that $\eps(l)<\eps(l{-}1)$, 
\begin{equation}\label{bedin}
\begin{split}
& \frac{1}{\eps^2}\int_{\Omega_1} 
\eps^{\alpha_H}H_*(\nabla_\eps^2\wh y_{l\eps})\dd x < \frac{1}{4l},
\quad 
\|\nabla_\eps\wh y_{l\eps}-(\nabla' \wh y_l|\nu_{\wh y_l})\|_{L^2}
< \frac{1}{4l}
\quad\text{ for all }\eps<\eps(l),\\
&
\Big|\frac{1}{\eps^2}\int_{\Omega_1} 
    W_\text{el}( x',\nabla_\eps  \wh y_{l\eps} ( x)  M_\eps( x)^{-1})\dd x
    -\frac{1}{2}\int_{\Omega_1}\overline{Q}_2(x',\nabla \wh y_l^\top\nabla \nu_{\wh y_l})\,\dd x'\Big|<\frac{1}{4l}\quad\text{ for all }\eps<\eps(l),
\end{split}
\end{equation}
and $\eps_l\downarrow 0$ as $l\to\infty$ which is possible 
due to Step~3 and Theorem~\ref{thm:dicht}.

For an arbitrary $\eps\in(0,\eps_0)$, we define 
\begin{equation}\label{approx}
\wh y_{\eps}:=\wh y_{l\eps}\quad\text{ if }\eps\in(\eps_{l+1},\eps_l]
\quad\text{
for this }l\in\N.
\end{equation} 
For $\eps\in(\eps_{l+1},\eps_l]$,  Step~1 and \eqref{bedin} ensure that for $l\in \mathbb{N}$
\[
\Big|
\frac{1}{\eps^2}\int_{\Omega_1}\Big\{ 
    W_\text{el}( x',\nabla_\eps  \wh y_{\eps} ( x)  M_\eps( x)^{-1})
    +\eps^{\alpha_H}H_*(\nabla_\eps^2\wh y_{\eps})\Big\}
    \dd x
    - \calM_0(\wh y)\Big| <\frac{1}{l}.
\]
Thus, letting $l\to \infty$ (implying that 
$(\eps_{l+1},\eps_l]\ni\eps\downarrow 0$), we find 
\begin{align*}
    &\frac{1}{\eps^2}\int_{\Omega_1}\Big\{ 
    W_\text{el}( x',\nabla_\eps  \wh y_{\eps} ( x)  M_\eps( x)^{-1})
    +\eps^{\alpha_H}H_*(\nabla_\eps^2\wh y_{\eps})\Big\}
    \dd x
    \to \calM_0(\wh y),\\
    &\text{and }
    \calM_\eps(\wh y_\eps)\to\calM_0(\wh y)
    \text{ as } \eps\to 0.
\end{align*}

Therefore, for suitable $c>0$ and $\wh\eps_0>0$ it is ensured that 
 $\calM_{\eps}(\wh y_{\eps})\le c$  for $\eps\in (0,\wh\eps_0]$. 
 For each such $\eps$, the function $\wh y_\eps\in W^{2,q_H}(\Omega_1,\R^3)$ has
 bounded elastic energy, therefore, the Healey--Kr\"omer 
 Theorem~\ref{thm:Healey} ensures that 
 $\det F_{\eps}>0$ a.e.\ in $\Omega_1$.
 Moreover, $\nabla_{\eps}\wh y_{\eps}\to (\nabla'\wh y|\nu_{\wh y})=:R_{\wh y}$ in $L^2(\Omega_1,\R^{3\times 3})$ with $R_{\wh y}\in \mathrm{SO(3)}$ a.e.\ in $\Omega_1$. 
This is guaranteed by the fact that for $\eps\in(\eps_{l+1},\eps_l]$,
we have by Step 1, \eqref{bedin}, and \eqref{approx}
\[
\|\nabla_\eps \wh y_\eps- (\nabla' \wh y|\nu_{\wh y})\|_{L^2}
\le 
\big\|\nabla_\eps \wh y_\eps- (\nabla' \wh y_l|\nu_{\wh y_l})\big\|_{L^2}
+
\big\|(\nabla' \wh y_l|\nu_{\wh y_l})-(\nabla' \wh y|\nu_{\wh y})\big\|_{L^2}
\le \frac{1}{l},\quad l\in\N. 
\]
Moreover, $\wh y_{\eps}\to \wh y$ in $W^{1,2}(\Omega_1,\R^3)$ holds.

5. For the electrical part of the energy, we proceed 
as follows.
Since the permittivity tensor $\bfk$ is positive definite, we have the Cholesky decomposition 
$\bfk(x)=L(x)L(x)^\top$ with a regular lower triangular matrix $L(x)\in\R^{3\times 3}$.
Similarly to \eqref{eq:hilfa} in Proposition~\ref{prop:conv}, but now with 
\begin{equation}\label{eq:PsiT}
\Psi:\Omega_1\times\R^{3\times 3}\to \R^{3\times 3},\quad
\Psi(x,F):=\begin{cases}
L(x)^\top F^{-\top}\sqrt{\det F}& \text {if }\mathrm{dist}(F,\SO(3))\le \frac{1}{2},\\
0 & \text{ else},
\end{cases}
\end{equation} 
and for
$Z=(\nabla'\varphi_0,m_0)$, 
$Z_{\eps}=\nabla_{\eps}\varphi_{\eps}$, and $\wh F_{\eps}=\nabla_{\eps} \wh y_{\eps}$, we have  
$Z_{\eps} \rightharpoonup Z$ in $L^{p_W}(\Omega_1,\R^3)$ by assumption and
$\wh F_{\eps}\to R_{\wh y}$ in $L^2(\Omega_1,\R^{3\times 3})$, where 
$\det \wh F_{\eps}>0$ and $R_{\wh y}\in \SO(3)$ a.e.\ in $\Omega_1$.
We derive the estimate
\begin{multline*}
\int_{\Omega_1}\bfk (R_{\wh y} Z)\cdot (R_{\wh y}Z)\,\dd x
=\int_{\Omega_1}(L^\top R_{\wh y} Z)\cdot (L^\top R_{\wh y}Z)\,\dd x\\
\le 
\liminf_{\eps\to 0}
\int_{\Omega_1} \bfk (\wh F_{\eps}^{-\top}Z_{\eps})\cdot (\wh F_{\eps}^{-\top}Z_{\eps})\det \wh F_{\eps}\,\dd x.
\end{multline*}
Therefore, we obtain by the definition of $m=m_{\wh y,\varphi_0}$ in  
$\calE_0(y_0,\varphi_0)$ and $\varphi_{\eps}\to \varphi_0$ in $L^2(\Omega_1)$,
\begin{equation*}
\begin{split}
\liminf_{{\eps}\to 0} \calE_{\eps}(\wh y_{\eps},\varphi_{\eps})
& =\frac{\beta}{2}\liminf_{{\eps}\to 0}
\int_{\Omega_1} (\bfk \wh F_{\eps}^{-\top}Z_{\eps})\cdot (\wh F_{\eps}^{-\top}Z_{\eps})\det \wh F_{\eps}\,\dd x
-\gamma \lim_{\eps\to 0} \int_{\Omega_1}n_\mathrm{ch}\varphi_\eps\,\dd x\\
&\ge \frac{\beta}{2}\int_{\Omega_1}R_{\wh y}^\top\bfk R_{\wh y}(\nabla'\varphi_0,m_0)
\cdot  (\nabla'\varphi_0,m_0)\,\dd x
-\gamma \int_{\Omega_1}n_\mathrm{ch}\varphi_0\,\dd x\\
&\ge \frac{\beta}{2}\int_{\Omega_1}R_{\wh y}^\top\bfk R_{\wh y}(\nabla'\varphi_0,m)
\cdot (\nabla'\varphi_0,m)\,\dd x
-\gamma \int_{\Omega_1}n_\mathrm{ch}\varphi_0\,\dd x
=\calE_0(y_0,\varphi_0).
\end{split}
\end{equation*}	

6. Finally, combining Steps 4 and 5, we verify 
\begin{equation*}
\begin{split}
\limsup_{\eps\to 0}\calF_{\eps}(\wh y_{\eps},\varphi_{\eps})
&\le\limsup_{\eps\to 0}\calM_{\eps}(\wh y_{\eps})
-\liminf_{\eps\to 0}\calE_{\eps}(\wh y_{\eps},\varphi_{\eps})
\le \calM_0(\wh y) - \calE_0(\wh y,\varphi_0)
=\calF_0 (\wh y,\varphi_0),
\end{split}
\end{equation*}
which completes the proof.
\end{proof}

Finally, we show that property (IV) holds.

\begin{lemma}\label{lem:III}{\rm{[Property (IV)].}}
We assume \ref{assu:W} -- \ref{eq:AssuElectric2}.
Let $(y_\eps,\varphi_\eps)_{\eps>0}\subseteq \calY_1 \times \calV_1$ 
be a sequence fulfilling 
 \eqref{eq:le}. Moreover,
let $y_0\in \calY_0$ and  
$\varphi_0\in\calV_0$ be the limits 
with respect to the convergences stated in Proposition~\ref{prop:conv}, 
 then 
\[
\liminf_{\eps\to 0} \calE_\eps(y_\eps,\varphi_\eps)\geq \calE_0(y_0,\varphi_0).
\]
\end{lemma}

\begin{proof}
Similar to Step~5 of the proof in Lemma~\ref{lem:II}  with the definition of $\Psi$ as in \eqref{eq:PsiT}
and 
$Z=(\nabla'\varphi_0,m_0)$, 
$Z_\eps=\nabla_\eps\varphi_\eps$, and $F_\eps=\nabla_\eps y_\eps,$  one derives the estimate
\begin{multline*}
\int_{\Omega_1}\bfk (R_{y_0} Z)\cdot (R_{y_0}Z)\,\dd x
=\int_{\Omega_1}(L^\top{R}_{y_0}Z)\cdot (L^\top R_{y_0}Z)\,\dd x
\\\le 
\liminf_{\eps\to 0}  
\int_{\Omega_1} \bfk (F_\eps^{-\top}Z_\eps)\cdot (F_\eps^{-\top}Z_\eps)\det F_\eps\,\dd x.
\end{multline*}
Therefore, we obtain by the definition of $m = m_{y_0,\varphi_0}$ in the limit energy 
$\calE_0(y_0,\varphi_0)$
\begin{equation*}
\begin{split}
\liminf_{\eps\to 0} \calE_\eps(y_\eps,\varphi_\eps)
& =\frac{\beta}{2}\liminf_{\eps\to 0}
\int_{\Omega_1} (\bfk F_\eps^{-\top}Z_\eps)\cdot (F_\eps^{-\top}Z_\eps)\det F_\eps\,\dd x
-\gamma \lim_{\eps\to 0} \int_{\Omega_1}n_\mathrm{ch}\varphi_\eps\,\dd x\\
&\ge \frac{\beta}{2}\int_{\Omega_1}{R_{y_0}^\top}\bfk R_{y_0}(\nabla'\varphi_0,m_0)
\cdot  (\nabla'\varphi_0,m_0)\,\dd x
-\gamma \int_{\Omega_1}n_\mathrm{ch}\varphi_0\,\dd x\\
&\ge \frac{\beta}{2}\int_{\Omega_1}R_{y_0}^\top\bfk R_{y_0}(\nabla'\varphi_0,m)
\cdot (\nabla'\varphi_0,m)\,\dd x
-\gamma \int_{\Omega_1}n_\mathrm{ch}\varphi_0\,\dd x
=\calE_0(y_0,\varphi_0),
\end{split}
\end{equation*}								
which is the desired result.
\end{proof}

\subsection{Some convergence results for subproblems}\label{ss34}

Here we collect convergence properties that are relevant in finding recovery sequences as stated in Lemma~\ref{lem:IV} and Lemma~\ref{lem:II}, 
respectively.

\begin{lemma}[Convergence of pulled-back tensor]
\label{lem:qw6}
We assume that \ref{assu:W} -- \ref{eq:AssuElectric2} are satisfied.
Let $(y_\eps)_{\eps>0}\subseteq \calY_1 $ 
be a sequence and  $y\in \calY_0$ and denote $R:=(\nabla'y|\nu)$ and $F_\eps=\nabla_\eps y_\eps$.
Let $F_\eps\to R$ in $L^2(\Omega_1,\R^{3\times 3})$ with $R\in \SO(3)$ a.e.\ in $\Omega_1$.
Moreover, for some $q_W>6$, let
$\|F_\eps\|_{L^{q_W}}$, $\|R\|_{L^{q_W}}\le c$ and 
$\|(\det F_\eps)^{-1}\|_{L^{q_W/2}}\le c$, then 
\[
 F^{-1}_\eps\bfk F_\eps^{-\top}\det F_\eps
=\frac{\Cof F_\eps^\top\bfk \Cof F_\eps}{\det F_\eps }
\to 
R^\top\bfk \, R\text{ in } 
L^1(\Omega_1,\R^{3\times 3}).
\] 
\end{lemma}

\begin{proof}
1. Fix $\theta\in(0,1)$. From $\|F_\eps- R\|_{L^q}\le \|F_\eps- R\|_{L^{q_W}}^{1-\theta}
\|F_\eps- R\|_{L^2}^{\theta}$ for $1/q=(1{-}\theta)/q_W+\theta/2$, we obtain also
$F_\eps\to R$ in $L^q(\Omega_1,\R^{3\times 3})$ for all $q\in[2,q_W)$.

In three dimensions, the determinant of a matrix can be expressed as a sum of terms, each being a product of three matrix entries. Consequently, 
$
\det F_\varepsilon \to \det R \quad \text{in } L^{q/3}(\Omega_1)$ holds for all  $q \in [3, q_W)$.
Similarly, the cofactor matrix has entries that are sums of products of two matrix entries, which implies the convergence
$\Cof F_\varepsilon \to \Cof R \quad \text{in } L^{q/2}(\Omega_1; \mathbb{R}^{3\times 3})$ for all 
$q \in [2, q_W)$.

2. Since $R\in \SO(3)$ a.e.\ in $\Omega_1$, 
we have $1/\det F_\eps-1/\det R=(1{-}\det F_\eps)/\det F_\eps$, and
\begin{equation*}
\begin{split}
\int_{\Omega_1}\Big|\frac{1}{\det F_\eps}-\frac{1}{\det R}\Big|\dd x
& 
=\int_{\Omega_1}\Big|\frac{1-\det F_\eps}{(\det F_\eps )}\Big|\dd x
\le \Big\|\frac{1}{(\det F_\eps )}\Big\|_{L^{\frac{q_W}{2}}}
\|1{-}\det F_\eps\|_{L^{\frac{q_W}{q_W-2}}}.
\end{split}
\end{equation*}
By assumption, we have $q_W>6$ such that $\frac{q_W}{q_W-2}<\frac{q_W}{3}$ (here we need in fact only $q_W>5$), 
and thus by Step 1,
$\det F_\eps\to 1$ in $L^{\frac{q_W}{q_W-2}}(\Omega_1)$.
Since we also assume that $\|\frac{1}{\det F_\eps}\|_{L^{q_W/2}}$ is bounded, it follows that
$
\frac{1}{\det F_\eps}\to 1$ in
$L^1(\Omega_1)$. We even have strong convergence in $L^t(\Omega_1)$ for all $t\in[1,q_W/2)$. Indeed, let again $\theta\in (0,1)$. From $\|z\|_{L^t}\le \|z\|_{L^{q_W/2}}^{1-\theta}
\|z\|_{L^1}^{\theta}$ for 
$1/t=(1{-}\theta)2/q_W+\theta$, we obtain
\[
\frac{1}{\det F_\eps}\to 1\text{  in }L^t(\Omega_1)
\text{ for all }t\in[1,q_W/2).
\]

3. For $q\in [6,q_W)$, it holds that $\frac{q}{q-4}\le \frac{q}{2}
< \frac{q_W}{2}$, and we estimate 
\begin{equation*}
\begin{split}
&\int_{\Omega_1}
\Big|\frac{\Cof F_\eps^\top\bfk \Cof F_\eps}{\det F_\eps }
- \Cof R^\top \bfk \Cof R\Big|\dd x\\
&
=\int_{\Omega_1}\Big|
\Big(\frac{1}{\det F_\eps}{-}1\Big)\Cof F_\eps^\top\bfk\Cof F_\eps
+(\Cof F_\eps^\top{-}\Cof R^\top)\bfk\Cof F_\eps\\
& \hspace{1.5cm}
+\Cof R^T\bfk(\Cof F_\eps{-}\Cof R)\Big|\dd x\\
& \le 
\Big\|\frac{1}{\det F_\eps}{-}1\Big\|_{L^\frac{q}{q-4}}
\|\Cof F_\eps\|_{L^{q/2}}
\|\bfk\|_{L^\infty}
\|\Cof F_\eps\|_{L^{q/2}}\\
& \hspace{1.5cm} +
\|1\|_{L^\frac{q}{q-4}}
\|\Cof F_\eps-\Cof R\|_{L^{q/2}}
\|\bfk\|_{L^\infty}
\|\Cof F_\eps\|_{L^{q/2}}\\
& \hspace{1.5cm}+
\|1\|_{L^\frac{q}{q-4}}
\|\Cof R\|_{L^{q/2}}
\|\bfk\|_{L^\infty}
\|\Cof F_\eps-\Cof R\|_{L^{q/2}}.
\end{split}
\end{equation*}
The convergences proven in Steps 1 and 2 and 
$\Cof R^\top \bfk \Cof R=
R^{-1}\bfk  R^{-\top}=R^\top\bfk  R$ for $R\in  \SO(3)$ ensure 
the $L^1$-convergence stated in the lemma.
\end{proof}

The proof of the next result is mainly inspired by ideas in \cite[Section 2]{ag-ago19}.

\begin{theorem}[Recovery sequence for the elastic energy]
\label{thm:dicht}
In addition to the Assumptions \ref{assu:W}-- \ref{assu:Prestrain}, we  suppose that the Assumption \ref{assu:B} holds. Then, for every $\wh y\in W^{2,2}_\mathrm{iso}(\omega,\R^3)
\cap C^\infty(\overline{\omega},\R^3)$ there exists a sequence of deformations $y_\eps\in W^{2,q_H}(\Omega_1,\R^3)$ such that 
\[
y_\eps\to \wh y\text{ in }W^{1,2}(\Omega_1,\R^3),\quad
\nabla_\eps y_\eps\to (\nabla'\wh y|\nu_{\wh y}) \text{ in } 
L^2(\Omega_1,\R^{3\times 3}),\text{ and }
\]
\[\lim_{\eps\to 0}\frac{1}{\eps^2}\calM_\eps^{\mathrm{el}}(y_\eps)
= \calM_0(\wh y), \quad\text{where }
\calM_\eps^{\mathrm{el}}(y_\eps):=
\int_{\Omega_1}W_\mathrm{el}(x', \nabla_\eps y_\eps M_\eps(x)^{-1})\dd x.
\]
Moreover, $\eps\|y_\eps\|_{W^{2,q_H}}(\Omega_1,\R^3)\to 0$ for $\eps\to 0$.
\end{theorem}

\begin{proof}
1. The Assumption \ref{assu:B} is used in \cite{ag-ago19} to derive models for heterogeneous elastic plates with in-plane modulation of the target curvature. 
In \cite [Lemma 2.2]{ag-ago19} it is proved that
\[
\argmin_{s\in \mathrm{Sym}(2)}\int_{-1/2}^{1/2}
Q_2\big(x',s+x_3G-B_{2\times 2}(x',x_3)\big)\dd x_3= \ol B_{2\times 2}(x'):=\int_{-1/2}^{1/2} 
B_{2\times 2}(x',t)\dd t
\]
for a.a.\ $x'\in\omega$ and all $G\in \mathrm{Sym}(2)$. In particular, the minimizer is independent of $G$. 

Let now the vector potential $g\in W^{1,2}(\omega;\R^2)$ be as in \eqref{vektorpotential}, i.e.\ $\nabla_\mathrm{sym} g = \ol B_{2\times 2}$, which of course also does not depend on $G$.
Therefore, we obtain
\begin{equation}\label{optc}
\begin{split}
\overline{Q}_2(x',G)
& =\int_{-1/2}^{1/2} Q_2\Big(x',x_3G
+\ol B_{2\times 2}(x')
-B_{2\times 2}(x',x_3)\Big)\dd x_3\\
& =\int_{-1/2}^{1/2} Q_2\Big(x',
x_3G+\nabla_{\mathrm{sym}}g-B_{2\times 2}(x',x_3)\Big)\dd x_3.
\end{split}
\end{equation}
Using e.g.\ Lemma~\ref{mollification} ii), we can approximate the potential $g\in W^{1,2}(\omega,\R^2)$ by $g_\eps\in W^{2,q_H}(\omega,\R^2)$ (also not depending on $G$) such that 
$\eps\|g_\eps\|_{W^{2,q_H}(\omega,\R^2)}\to 0$ for $\eps\to 0$. 

2. Let $\wh y\in W^{2,2}_\mathrm{iso}(\omega,\R^3)\cap C^\infty(\overline{\omega},\R^3)$ and the corresponding normal vector $\nu:=\nu_{\wh y}\in C^\infty(\overline{\omega},\R^3)$ be given. Then, we have
\[
R:=R_{\wh y}=(\nabla'\wh y|\nu_{\wh y})=(\nabla'\wh y|\nu)\in C^\infty(\overline{\omega},\R^{3\times 3})\text{ and }R\in \mathrm{SO}(3) \text{ a.e.\ in }\omega.
\]
For $d_\eps\in W^{2,q_H}(\Omega_1,\R^3)$, 
which also provide an approximation for a  $d\in L^2(\Omega_1,\R^3)$ that is determined below 
(cf.\ also Lemma~\ref{mollification} i)) 
and $g_\eps:=(g_{\eps 1},g_{\eps 2})\in W^{2,q_H}(\omega,\R^2)$ as above,  
we define the sequence 
\begin{equation}\label{ansatz}
y_\eps(x',x_3)=\wh y(x') + \eps[ x_3 \nu(x')+ \nabla' \wh y(x') g_\eps(x')] + \eps^2 
D_\eps(x',x_3),
\end{equation}
 where 
 \[
D_\eps(x',x_3)=\int_0^{x_3}d_\eps(x',t)\dd t.
\]
Note that the properties $\wh y\in W^{2,2}_\mathrm{iso}(\omega,\R^3)\cap C^\infty(\overline{\omega},\R^3)$ 
(yielding in particular $|\partial_1 \wh y|=|\partial_2 \wh y|=1$),
$\nu \in C^\infty(\overline{\omega},\R^3)$, and 
$D_\eps \in W^{2,q_H}(\Omega_1,\R^3)$ imply that
$y_\eps\in W^{2,q_H}(\Omega_1,\R^3)$. Moreover, we easily confirm that  
$y_\eps\to \wh y$ in $W^{1,2}(\Omega_1,\R^3)$ for $\eps \to 0$.  
Moreover, we have the explicit estimate
\begin{equation}\label{convdicht}
\eps \|y_\eps\|_{W^{2,q_H}}
\le \eps c \Big(\|\wh y\|_{C^\infty}
+\eps \big(\|\wh y\|_{C^\infty}
+\|\wh y\|_{C^\infty}\|g_\eps\|_{W^{2,q_H}}\big)
+\eps^2\|d_\eps\|_{W^{2,q_H}}\Big)\to 0,
\end{equation}
which ensures the last assertion of the theorem.

A direct computation shows that
\begin{equation}\label{epsgrad}
\begin{split}
\nabla_\eps y_\eps(x',x_3)&=R(x')
+\eps\Big( x_3 \nabla'\nu(x')+\nabla'[\nabla' \wh y(x')g_\eps(x')]\big|d_\eps(x)\Big)
+\eps^2
(\nabla'D_\eps(x)|0),\\
\nabla_\eps y_\eps(x',x_3)M_\eps(x)^{-1}
&=\Big(R(x')
+\eps\big( x_3 \nabla'\nu(x')+\nabla'[\nabla' \wh y(x')g_\eps(x')]\big|d_\eps(x)\big)
+\eps^2
(\nabla'D_\eps(x)|0) 
 \Big)\\&\hspace{2em}\times
(\mathbb{I}_3-\eps B(x)+\eps^2B(x)^2)
+O(\eps^3).
\end{split}
\end{equation}
Since $R\in\mathrm{SO(3)}$, the frame indifference of $W_\mathrm{el}$ yields 
\begin{equation*}
\begin{split}
&W_\mathrm{el}(x',\nabla_\eps y_\eps M_\eps^{-1})
=
W_\mathrm{el}(x',R^\top\nabla_\eps y_\eps(\cdot,x_3)M_\eps^{-1})\\
&=W_\mathrm{el}\Big(x',\Big\{\mathbb{I}_3
+\eps R^\top\big( x_3 \nabla'\nu+\nabla'[\nabla' \wh yg_\eps]\big|d_\eps(\cdot,x_3)\big)
+\eps^2 R^\top
(\nabla'D_\eps(\cdot,x_3)|0) 
 \Big\}\\&\quad\times
(\mathbb{I}_3-\eps B(\cdot,x_3)+\eps^2B(\cdot,x_3)^2)+O(\eps^3)\Big)
.
\end{split}
\end{equation*}
For $\varepsilon$ sufficiently small,
 the properties \eqref{dconv} of the approximation 
 $d_\eps$ and $g_\eps$ guarantee
\begin{equation*}
\begin{split}
& \mathrm{dist}\big(R^\top\nabla_\eps y_\eps(x',x_3)M_\eps^{-1},\mathrm{SO}(3)\big)
\\
& \le \eps \|R\|_{L^\infty}\Big(
\|\nabla'\nu\|_{L^\infty}
+\|\wh y\|_{W^{2,\infty}} \|g_\eps\|_{L^\infty}
+\|\wh y\|_{W^{1,\infty}} \|g_\eps\|_{W^{1,\infty}}
+\|d_\eps\|_{L^\infty}
+\|B\|_{L^\infty}\Big)\\
& \quad +\eps^2\|R\|_{L^\infty}\Big(
\|\nabla'\nu\|_{L^\infty}
+\|\wh y\|_{W^{2,\infty}} \|g_\eps\|_{L^\infty}
+\|\wh y\|_{W^{1,\infty}} \|g_\eps\|_{W^{1,\infty}}
+\|d_\eps\|_{L^\infty}\Big)\|B\|_{L^\infty}\\
& \quad 
+\eps^2 \|R\|_{L^\infty}\|\nabla'D_\eps\|_{L^\infty}
+\eps^2\|B\|_{L^\infty}^2
+K\eps^3<\epsilon_0, 
\end{split}
\end{equation*}
where $\epsilon_0>0$ is such that in the $\epsilon_0$-neighbourhood of $\mathrm{SO}(3)$, the elastic energy density $W_\mathrm{el}$ is bounded and $C^2$ regular.
Using a Taylor expansion of $W_\mathrm{el}$ at 
$\mathbb{I}_3$ we get
\begin{equation*}
\begin{split}
\frac{1}{\eps^2}\calM_\eps^{\mathrm{el}}(y_\eps)
&=\frac{1}{\eps^2}\int_{\Omega_1}
W_\mathrm{el}\Big(\mathbb{I}_3
+\eps R^\top\Big( x_3 \nabla'\nu+\nabla'[(\nabla' \wh y)g_\eps]\big|d_\eps\Big)
+\eps^2 R^\top
(\nabla'D_\eps|0) 
 \Big)\times\\
 &\hspace{2cm}\times
(\mathbb{I}_3-\eps B+\eps^2B^2)+O(\eps^3)
\Big)
\dd x\\
&=\frac{1}{2}\int_{\Omega_1}\Big\{
D^2W_\mathrm{el}(\mathbb{I}_3)
\Big(R^\top\big( x_3 \nabla'\nu+\nabla'[(\nabla' \wh y)g_\eps]\big|d_\eps(x)\big)-B\Big)^{\otimes 2}
+O(\eps) \Big\}\dd x,
\end{split}
\end{equation*}
resulting in
\begin{equation*}
\begin{split}
\lim_{\eps\to 0}
\frac{1}{\eps^2}\calM_\eps^{\mathrm{el}}(y_\eps)
&=\lim_{\eps\to 0}
\frac{1}{2}\int_{\Omega_1}
D^2W_\mathrm{el}(\mathbb{I}_3)
\Big(R^\top\big( x_3 \nabla'\nu+\nabla'[(\nabla' \wh y)g_\eps]\big|d_\eps\big)-B\Big)^{\otimes 2}
 \dd x
 \\
&=\lim_{\eps\to 0}
\frac{1}{2}\int_{\Omega_1}
Q_3\big(R^\top\big( x_3 \nabla'\nu+\nabla'[(\nabla' \wh y)g_\eps]\big|d_\eps\big)-B\big)
\dd x\\
&=\frac{1}{2}\int_{\omega}\int_{-1/2}^{1/2}
Q_3\big(
R^\top\big( x_3 \nabla'\nu+\nabla'[(\nabla' \wh y)g]\big|d(x)\big)-B
\big)
 \dd x_3 \dd x'.
\end{split}
\end{equation*}
Here we used $\wh y,\,\nu\in C^\infty(\overline{\omega};\R^3)$, 
 $d_\eps\to d$ in $L^2(\Omega_1,\R^3)$, $g_\eps\to g$ in $W^{1,2}(\omega,
 \R^2)$.

Now, we follow the ideas in the proof of \cite[Theorem 2.6 (III)]{ag-ago19}.
Since $Q_3(F)=Q_3(F_\mathrm{sym})$ for all $F\in\R^{3\times 3}$, we evaluate the symmetric part of the argument in $Q_3$, and we will write it in a special form 
(see \eqref{help1} below). Indeed, let $X$ denote the $2\times 2$ 
upper left part of the $3\times 3$ matrix in the argument in $Q_3$, and let $X^\circ$ be defined
as in \eqref{eq:circNotation}.
We consider the map $L:\mathrm{Sym}(2)\to \R^3$ that assigns to $X\in \mathrm{Sym}(2)$ the unique  vector  
$c_X:=\mathrm{argmin}_{c\in \R^3}Q_3\big(X^\circ+(c\otimes e_3)\big)$. 
Since $Q_3$ is quadratic, the first order optimality conditions 
for $c_X$ ensure the linearity of the mapping $L$. We introduce 
\[
\bar d(x):=R(x')\Big(L(X_\mathrm{sym}(x))-\Big(
\begin{array}{c}
(\nabla'(\nabla'\wh y(x)g(x')))^\top\nu(x')\\
0\\
\end{array}\Big) 
+(2B_{13}\,2B_{23}\,B_{33})(x)^\top\Big).
\]
We recall that due to $X\in L^\infty(\Omega_1,\mathrm{Sym}(2))$, 
$B\in L^\infty(\Omega_1,\mathrm{Sym}(3))$, $\wh y,\, \nu\in C^\infty(\bar\omega;\R^3)$, and $g\in W^{1,2}(\omega,\R^2)$, thus,
$\bar d\in L^\infty(\Omega_1,\R^3)$ holds.
Using the vector $d=\bar d$  in the respective position 
in $Q_3$, 
we establish a representation of the symmetric part of the argument in $Q_3$ by
\begin{equation}\label{help1}
\Big(R^\top\big( x_3 \nabla'\nu(x')+\nabla'[\nabla' \wh y(x')g(x')]\big|\bar d(x)\big)-B\Big)_\mathrm{sym}
=\left(\begin{array}{cc}
&0\\[-1.5ex]
X_\mathrm{sym}&\\[-1.5ex]
&0\\[-0.5ex]
0\quad 0 & 0\\
\end{array}\right)
+(L(X_\mathrm{sym}(x)) \otimes e_3)_\mathrm{sym}.
\end{equation}
Due to the identity
$(\nabla' \wh y)g=g_1\partial_1 \wh y + g_2\partial_2 \wh y$, we find
\[
\nabla'\big((\nabla' \wh y)g\big)=
\big(g_1\partial_1 \partial_1 \wh y
+g_2\partial_1 \partial_2 \wh y\big|
g_1\partial_2 \partial_1 \wh y
+g_2\partial_2 \partial_2 \wh y\big)
+\nabla' \wh y\nabla' g.
\]
For $\wh y\in W^{2,2}_{\mathrm{iso}}(\omega;\R^3)$, 
we have 
$\partial_i \wh y\cdot \partial_j \wh y=\delta_{ij}$ and 
$\partial_i  \partial_j\wh y\cdot \partial_k\wh y=0$ for $i,j,k=1,2$.
Moreover, $\sum_{j=1}^3\nu_j\partial_i \nu_j=\tfrac{1}{2}\partial_i|\nu|^2=0$, $i=1,2$, since $\nu$ has norm one. These identities lead to
\[
(\nabla' \wh y)^\top \nabla' \big((\nabla' \wh y)g\big)=\nabla' g,
\quad 
R^\top(\nabla' \nu|0)
 =\left(\begin{array}{cc}
 (\nabla' \wh y)^\top \nabla' \nu &0\\
 0&0\\
 \end{array}\right).
\]
Therefore, we calculate 
\begin{equation}\label{help2}
X_\mathrm{sym}(x)=x_3(\nabla' \wh y(x'))^\top \nabla'\nu(x')
+\nabla'_{\mathrm{sym}} g(x')-B(x)_{2\times 2}\quad\text{ a.e.\ in }
\Omega_1.
\end{equation}
 Using \eqref{help1}, \eqref{help2},  and the definitions of $Q_2$ and $\overline{Q}_2,$
 we establish 
\begin{equation*}
\begin{split}
\lim_{\eps\to 0}
\frac{1}{\eps^2}
\calM_\eps^{\mathrm{el}}(y_\eps)& =
\lim_{\eps\to 0}
\frac{1}{\eps^2}\int_{\Omega_1} W_\mathrm{el}(x', \nabla_\eps y_\eps M_\eps(x)^{-1})\dd x\\ 
& =
\frac{1}{2}\int_{\omega}\int_{-1/2}^{1/2}
Q_3\big(x',
R^\top\big( x_3 \nabla'\nu(x')+\nabla'[\nabla' \wh y(x')g(x')]\big|\bar d(x)\big)-B(x)
\big)
 \dd x_3 \dd x'\\
 & =
\frac{1}{2}\int_{\omega}\int_{-1/2}^{1/2}
Q_2\big(x',
x_3\nabla' \wh y^\top \nabla' \nu+\nabla'_{\mathrm{sym}}g -B(x)_{2\times 2}
\big)
 \dd x_3 \dd x'\\
& = \frac{1}{2}\int_{\omega}\min_{s\in\R^{2\times 2}}\int_{-1/2}^{1/2}
Q_2\big(x',
x_3\nabla' \wh y^\top \nabla' \nu+s-B(x)_{2\times 2})
 \dd x_3 \dd x'\\
&= \frac{1}{2}\int_{\omega}
\overline{Q}_2(x',
\nabla \wh y^\top \nabla \nu)
 \dd x'
 =\calM_0(\wh y).
 \end{split}
 \end{equation*}
Finally, \eqref{epsgrad}, $\eps\|d_\eps\|_{W^{2,q_H}(\Omega_1,\R^3)}\to 0$,
and $\eps\|g_\eps\|_{W^{2,q_H}(\omega,\R^2)}\to 0$ for $\eps\to 0,$ ensure
that $\nabla_\eps y_\eps\to (\nabla'\wh y|\nu_{\wh y})$ in 
$L^2(\Omega_1,\R^{3\times 3})$, 
which completes the proof.
 \end{proof}

\section{Concluding remarks}\label{sec:conclusion}
In this manuscript, we have restricted our considerations to prestrains satisfying \ref{assu:B}, primarily for the sake of clarity and transparency in the proof of Theorem~\ref{thm:dicht}. 
We note, however, that the work~\cite{Padilla-Garza2022} addresses more general prestrains that satisfy only \ref{assu:Prestrain}.

We expect that, also in this more general setting, and in the spirit of Theorem~\ref{thm:dicht}, 
for all deformations  $\wh y\in W^{2,2}_{\mathrm{iso}}(\omega,\R^3)
\cap C^\infty(\overline{\omega},\R^3)$ a sequence 
$(y_\eps)_{\eps>0}$ can be constructed
with $y_\eps\in W^{2,q_H}(\Omega_1,\R^3)$,
$
y_\eps\to \wh y$ in $W^{1,2}(\Omega_1,\R^3)$,
$\nabla_\eps y_\eps\to (\nabla'\wh y|\nu_{\wh y})$ in  
$L^2(\Omega_1,\R^{3\times 3})$, and
$\lim_{\eps\to 0}\frac{1}{\eps^2}\calM_\eps^{\mathrm{el}}(y_\eps)
= \calM_0(\wh y)$, and the additional property that 
$\eps\|y_\eps\|_{W^{2,q_H}(\Omega_1,\R^3)}\to 0$ for $\eps\to 0$ 
is satisfied. This would then ensure the validity of 
Lemma~\ref{lem:II} and Theorem~\ref{thm:main}
also in the more general prestrain setting.

Following the strategy of the proof of~\cite[Theorem~3]{Padilla-Garza2022}, 
one would have to distinguish three cases for the second fundamental form 
$\widehat{\mathrm{II}} := (\nabla \widehat{y}^\top \nabla \nu_{\widehat{y}})$, namely:  
\begin{enumerate}[label=(\roman*)]
\item $\widehat{\mathrm{II}} = 0$, treated by means of the ansatz in Eqn.\ (93) in~\cite{Padilla-Garza2022},  
\item $\widehat{\mathrm{II}}$ bounded away from zero, corresponding to the ansatz in Eqn.~(115) in~\cite{Padilla-Garza2022}, and  
\item the general case, where $\widehat{\mathrm{II}}$ is neither identically zero nor nowhere vanishing.  
\end{enumerate}
In the latter case, the previous two ansatzes would need to be smoothly connected 
by means of a transition layer along the boundary of the set $\{\widehat{\mathrm{II}} = 0\}$. 
This could be achieved through an ansatz of the intricate form in Eqn.\ (139) in~\cite{Padilla-Garza2022}.

In addition, one must ensure that all terms in the sum belong to the space
$W^{2,q_H}(\Omega_1, \mathbb{R}^3)$ and that $
\varepsilon \|y_\varepsilon\|_{W^{2,q_H}(\Omega_1, \mathbb{R}^3)} \to 0$ as $\eps\to0$.
This, in particular, would require an additional smoothing procedure 
for the various factors appearing in the different terms of~(139) 
in~\cite{Padilla-Garza2022}. 
We avoided this technically involved construction by adopting the simplifying assumption~\ref{assu:B}.

Finally, let us remark that also in the limit model a result similar to Corollary \ref{cor:FE} holds.

\begin{lemma}
We assume \ref{assu:W} -- \ref{eq:AssuElectric2}.
Let $y_0\in \calY_0$ and  
$\varphi_0\in\calV_0$ be the limits 
with respect to the convergences stated in Proposition~\ref{prop:conv}.
If $\calF_0(y_0,\varphi_0)\ge \calF_0(y_0,\overline{\varphi})$ for all
$\overline{\varphi}\in\calV_0,$ meaning
$\calE_0(y_0,\varphi_0)\le \calE_0(y_0,\overline{\varphi})$ for all
$\overline{\varphi}\in\calV_0,$ then
\[
\calN_0(\varphi_0)=2\calQ_0(y_0,\varphi_0), \quad 
\calE_0(y_0,\varphi_0)=\calQ_0(y_0,\varphi_0)-\calN_0(\varphi_0),
\]
where 
\[
    \calQ_0(y_0,\varphi_0) 
    := \frac{\beta}{2}\int_{\omega} 
    \mathbb{K}^\mathrm{eff}_{y_0}(x')\nabla'\varphi_0( x')\cdot \nabla' \varphi_0( x')\dd  x',
    \quad
    \calN_0(\varphi_0):=\gamma\int_{\omega}  
    \ol n_\mathrm{ch}( x')  \varphi_0( x') \dd  x'.
\]
\end{lemma}
\begin{proof}
We use $\overline{\varphi}\in\calV_0$ of the form $\overline{\varphi}=\varphi_0(1+a)$, $a\in \R$. 
Then $\calE_0(y_0,\varphi_0)\le \calE_0(y_0,\overline{\varphi})$ for 
these $\overline{\varphi}$ yields
\[
\frac{\beta(2a+a^2)}{2}\int_\omega\mathbb{K}^\mathrm{eff}_{y_0}
\nabla'\varphi_0\cdot \nabla'\varphi_0\dd x'
-\gamma a \int_{\omega} \ol n_\mathrm{ch}\varphi_0\dd x'\le 0.
\]
Therefore, we find 
\begin{equation}\label{eq:1}
\gamma\int_\omega \ol n_\mathrm{ch}\varphi_0 \dd x'
 -\beta(1+\frac{a}{2})\int_\omega 
\mathbb{K}^\mathrm{eff}_{y_0}
\nabla'\varphi_0 \cdot \nabla'\varphi_0
\dd x'
\begin{cases}
\quad\ge 0 & \forall a<0,\\
\quad\le 0 & \forall a>0.
\end{cases}
 \end{equation}
We consider in \eqref{eq:1} the limit $a\uparrow 0$ and obtain 
\[
\gamma\int_\omega \ol n_\mathrm{ch}\varphi_0 \dd x'
 -\beta\int_\omega 
\mathbb{K}^\mathrm{eff}_{y_0}
\nabla'\varphi_0 \cdot \nabla'\varphi_0
\dd x' \ge 0.
 \]
For the opposite inequality we take the limit $a\downarrow 0$ in \eqref{eq:1}. 
This together ensures $\calN_0(\varphi_0)=2\calQ_0(y_0,\varphi_0)$. Relation \eqref{eq:limitElecEnergy} then finalizes the proof.
\end{proof}

\appendix

\section{Auxiliary tools}

The following result, due to Healey--Kr\"omer \cite{HeaKro2009IWSS} 
can also be found e.g.\ in\ \cite[Theorem 2.5.3]{ag-kru19} and \cite[Theorem 3.1]{ag-mie20}. 
\begin{theorem}[Healey--Krömer]\label{thm:Healey}
Assume that the elastic energy $W_\mathrm{el}:\Omega_1\times \R^{3\times 3}\to[0,\infty]$ satisfies
\ref{assu:Wgrowth}, and the hyperstress potential $H:\R^{3\times 3\times 3}\to \R_+$ fulfils \ref{assu:hyper} with $q_W \geq 3$ and $
q_W/2>3q_H/(q_H-3)$. 
Then, for all $C_\mathrm{M}>0$ there exists a $C_\mathrm{HK}=C_\mathrm{HK}(C_\mathrm{M},q_H,q_W)>0$ such that for all deformations $y\in W^{2,q_H}(\Omega_1,\R^3)$ with bounded elastic energy
$\int_{\Omega_1}\big(W_\mathrm{el}(\nabla y) + H(\nabla^2 y)\big)\dd x < C_\mathrm{M}$ it is satisfied
\begin{gather*}
\|y\|_{W^{2,q_H}}\leq C_\mathrm{HK},\quad \|y\|_{C^{1-3/q_H}}\leq C_\mathrm{HK},\quad \|(\nabla y)^{-1}\|_{C^{1-3/q_H}}\leq C_\mathrm{HK},\\
        \det\nabla y(x)\geq 1/C_\mathrm{HK}\text{ for all }x\in\Omega_1.
\end{gather*}
\end{theorem}

\begin{lemma}\label{mollification}
Let $q_H$ be defined as in Assumption \ref{assu:hyper}.

i) For  $d\in L^2(\Omega_1,\R^3)$, there exists a sequence $(d_\eps)_{\eps>0}$
such that $d_\eps\in W^{2,q_H}(\Omega_1,\R^3)$, $d_\eps\to d$ in $L^2(\Omega_1,\R^3)$, and 
$\eps\|d_\eps\|_{W^{2,q_H}}\to 0$ for $\eps\to 0$.

ii) For  $g\in W^{1,2}(\omega,\R^2)$, there exists a sequence $(g_\eps)_{\eps>0}$ such that $g_\eps\in W^{2,q_H}(\omega,\R^2)$,  
$g_\eps\to g$ in $W^{1,2}(\omega,\R^2)$, and
$\eps\|g_\eps\|_{W^{2,q_H}(\omega,\R^2)}\to 0$ for $\eps\to 0$.

\end{lemma}

\begin{proof} 
1. We fix $d\in L^2(\Omega_1,\R^3)$   
and construct approximations $d_{\eps}\in W^{2,q_H}(\Omega_1,\R^3)$ by solving 
the minimization problem
\[
\min_{\bar d\in W^{2,q_H}(\Omega_1,\R^3)}\calI_\eps(\bar d),\quad \text{where}\quad\calI_\eps(\bar d):=\int_{\Omega_1}\Big\{\frac{\eps^\tau}{q_H}\big\{|\nabla^2 \bar d|^{q_H}+|\nabla \bar d|^{q_H}\big\}
+\frac{1}{2}|\bar d{-}d|^2\Big\}\dd x
\]
for some $\tau\in(0,q_H)$. 
Note that the functionals $\calI_\eps$ 
are strictly convex and the minimizers, denoted by $d_\eps$, are unique.
The associated Euler--Lagrange equation reads
\[
0=\int_{\Omega_1}\Big\{\eps^{\tau}\big\{|\nabla^2 d_\eps|^{q_H-2}\nabla^2 d_\eps\,\,\vdots\,\nabla^2 \xi
+|\nabla d_\eps|^{q_H-2}\nabla d_\eps : \nabla \xi\big\}+ (d_\eps{-} d)\cdot\xi\Big\}
\dd x
\]
for all $\xi\in W^{2,q_H}(\Omega_1,\R^3)$.
Testing the latter by $\xi=d_{\eps}$, we find  
$\|d_{\eps}\|_{W^{2,q_H}(\Omega_1,\R^3)}^{q_H}\le c/\eps^{\tau}$. 
Since $\tau<q_H$, we therefore get
\begin{equation}\label{dconv}
\eps\|d_\eps\|_ {W^{2,q_H}}\le \wt c \eps^{1-\tau/q_H}\to 0
\text{ and }
\eps\|d_\eps\|_ {W^{1,\infty}}\to 0\text{ for }\eps\to 0.
\end{equation}
 Next, we show that
$d_\eps\rightharpoonup d$ in $L^2(\Omega_1,\R^3)$ by proving the $\Gamma$-convergence of the functional $\calI_\eps$ to $\calI_0$ 
in the weak $L^2$-topology, where the limiting functional $\calI_0$ is defined by  
$\bar d\mapsto \frac{1}{2}\int_{\Omega_1}|\bar d{-}d|^2\dd x$. Trivially, the unique minimizer 
of $\calI_0$ in $L^2(\Omega_1,\R^3)$ is given by $ d$.

First, it is easy to see that for all sequences 
$\bar d_\eps\rightharpoonup d^*$ in $L^2(\Omega_1,\R^3)$, we have 
$\liminf_{\eps\to 0}\calI_\eps(\bar d_\eps)
\ge \calI_0(d^*)$. Indeed, this estimate follows from $\calI_\eps(\bar d_\eps)\geq \calI_0 (\bar d_\eps)$ and the weak lower semi-continuity of the norm.

Second, for all $d^*\in L^2(\Omega_1,\R^3)$, we can  find a recovery sequence $\wh d_\eps\to d^*$ in $L^2(\Omega_1,\R^3)$ such that $\limsup_{\eps\to 0}\calI_\eps(\wh d_\eps)$ $
=\calI_0(d^*)$. Namely, if $d^*\in W^{2,q_H}(\Omega_1,\R^3)$, we take the constant sequence
$\wh d_\eps\equiv d^*,$ giving $\limsup_{\eps\to 0}\calI_\eps(d^*)=\calI_0(d^*)$. 
In case $d^*\notin W^{2,q_H}(\Omega_1,\R^3)$, by density arguments 
there exist $(d_k)\subset W^{2,q_H}(\Omega_1,\R^3)$ with $d_k\to d^*$ in $L^2(\Omega_1,\R^3)$ and $a_k:=\|d_k\|_{W^{2,q_H}}^{q_H},$ possibly tending to infinity. For all $\eps>0,$ there exists $k(\eps)\in\N$, $k(\eps)\ge\max_{\eps'>\eps} k(\eps')$ such that $a_{k(\eps)}\le \eps^{-\beta}$. Then, setting $\wh d_\eps:=d_{k(\eps)}$, since $\wh d_\eps=d_{k(\eps)}\to d^*$ in $L^2(\Omega_1,\R^3)$ for $\eps\to 0$, we find  
\begin{equation}
\begin{split}
\limsup_{\eps\to 0}\calI_\eps(\wh d_\eps)
&\le \limsup_{\eps\to 0}\Big\{
\eps^{1+\tau} a_{k(\eps)}+\|\wh d_\eps- d\|_{L^2}^2\Big\}\\
&\le \limsup_{\eps\to 0}\Big\{\eps+\|\wh d_\eps- d^* + d^* - d\|_{L^2}^2\Big\}
=\| d -d^*\|_{L^2}^2=\calI_0(d^*).
\end{split}
\end{equation}
The fundamental theorem of $\Gamma$-convergence gives the  weak convergence $d_\eps\rightharpoonup  d$ in $L^2(\Omega_1,\R^3)$. 

It remains to show that the convergence of $d_\eps$ is actually strong in $L^2(\Omega_1,\R^3)$.  Testing the Euler-Lagrange equation by $d_\eps$ gives the identity
\[
\|d_\eps\|_{L^2}^2 + \eps^\tau\big(\|\nabla^2 d_\eps\|_{L^{q_H}}^{q_H }
+
\|\nabla d_\eps\|_{L^{q_H}}^{q_H }
\big)= \int_{\Omega_1} d_\eps  d\dd x.
\]
Using the weak convergence of $d_\eps$, we can pass to the limit in the right-hand side, giving
\[
\| d\|_{L^2}^2
\leq \liminf_{\eps\to 0}\big\{
\|d_\eps\|_{L^2}^2 + \eps^\tau\big(\|\nabla^2 d_\eps\|_{L^{q_H}}^{q_H }
+
\|\nabla d_\eps\|_{L^{q_H}}^{q_H }
\big)
\big\}\leq \| d\|_{L^2}^2.
\]
We conclude that $d_\eps$ is strongly converging to $ d$. Moreover, we obtain the convergence
$\eps^\tau(\|\nabla^2 d_\eps\|_{L^{q_H}}^{q_H }
+
\|\nabla d_\eps\|_{L^{q_H}}^{q_H })\to 0$.

2. The approximation for $g\in W^{1,2}(\omega;\R^3)$ is obtained similar to Step 1 by solving 
\[
\min_{\bar g\in W^{2,q_H}(\omega,\R^2)}\calI_\eps^g(\bar g),\quad \calI_\eps^g(\bar g):=\int_{\omega}\Big\{\frac{\eps^\tau}{q_H}\big\{|\nabla^2 \bar g|^{q_H}+|\nabla \bar g|^{q_H}\big\}
+\frac{1}{2}|\nabla(\bar g- g)|^2
+\frac{1}{2}|\bar g- g|^2\Big\}\dd x',
\]
and proceeding as in Step 1.
\end{proof}

\subsection*{Acknowledgments}
Support of the Deutsche Forschungsgemeinschaft under Germany’s Excellence Strategy The Berlin
Mathematics Research Center MATH+ (EXC-2046/1,
project 390685689) is gratefully acknowledged.
The authors also thank David Padilla-Garza (Jerusalem) and Danka Lučić (Jyväskylä) for  fruitful discussions and helpful comments during the preparation of this manuscript.

\bibliography{Refs.bib}

\end{document}